\input amstex      
 \documentstyle{amsppt}
\magnification 1200 
   
\hoffset= -.05 truecm
\hsize= 16.5 truecm
\vsize= 21.95 truecm
\voffset= -.05 truecm

 \def\ns{\hskip -.15 truecm} 
 \NoBlackBoxes

\def\summary{Theorem }
\def\classthms{Theorems 3.4, 3.5, 4.1 and 4.2 }
\def\singlemmas{Theorems 2.5 and 2.6 }

 \def\notation{Notation }
 \def\defcan{Definition 1.1 } 
\def\hypsplit{\GPtrans \ns, Lemma 2.3 }   
\def\Galquad{Proposition 1.2 }
\def\Galquadone{(1.2.1) }
\def\biproconv{Proposition 1.3 }
\def\Galcycliconv{Proposition 1.4 }

\def\diagram{Definition 2.1 }
\def\defXbar{\diagram}
\def\diagramone{(2.1.1) }
\def\samegrplem{Lemma 2.2 }
\def\diagramGal{\samegrplem}
\def\ratlemma{Lemma 2.3 }
\def\discrlemma{Lemma 2.4 } 
 \def\prsinglemmaone{(2.4.1) }
 \def\prsinglemmatwo{(2.4.2) }
 \def\singlemma{Theorem 2.5 }
 \def\diagramtwo{(2.5.1) }
 \def\prsinglemmathree{(2.5.2) }
 \def\prsinglemmafour{(2.5.3) }
 \def\prsinglemmafive{(2.5.4) }
 \def\prsinglemmasix{(2.5.5) }
 \def\prsinglemmasixb{(2.5.6) }
 \def\thmsinglept{Theorem 2.6 }
\def\singleptone{(2.6.1) }
 \def\prsinglemmaseven{(2.6.2) }
 \def\prsinglemmaeight{(2.6.3) }
 \def\prsinglemmanine{(2.6.4) }

\def\convention{Definition 3.1 }
\def\convnumber{(3.1.1) }
\def\convnumb{\convnumber}
\def\split{Proposition 3.2 }

\def\biprocyc{Proposition 3.3 } 
\def\bipro{\biprocyc \ns, 1) }

 \def\crepantbid{Theorem 3.4 }

 \def\crepantcyc{Theorem 3.5 }

 \def\noncrepantbid{Theorem 4.1 }
 \def\ncrone{(4.1.1) }
 \def\ncrdiag{(4.1.2) }
\def\ncrbidthree{(4.1.3) }
\def\ncrfour{(4.1.4) }
\def\noncrepantcyc{Theorem 4.2 }
 \def\ncrtwo{(4.2.1) }
 \def\ncrthree{(4.2.2) }
\def\ncrfive{(4.2.3) }
\def\ncrsix{(4.2.4) }
 \def\ncrseven{(4.2.5) }
\def\classing{Theorems 3.4, 3.5, 4.1 and 4.2 }
\def\singsplit{Corollary 4.3 }

 \def\numbsing{Corollary 5.1 }
 \def\excrepantbid{Proposition 5.2 }
 \def\excrepantcyc{Proposition 5.3 }
 \def\exampnoncrepbid{Proposition 5.4 }
 
 \def\exampcyc{\GPsmooth \ns, Proposition 3.11 }

\def\Ca{[Ca] }
\def\CKM{[CKM] }
\def\GPtrans{[GP1] }
\def\GPsmooth{[GP2] }
\def\CR{[GP3] }
 \def\Hoone{[Ho1] }
 \def\Hothree{[Ho2] }
 \def\HM{[HM] }
 \def\Kon{[Ko] }
 \def\Pa{[Pa] }
 \def\Pu{[Pu] }

 \topmatter
 \title Classification of quadruple Galois canonical covers II
 \endtitle
 
 \author Francisco Javier Gallego \\ and \\ Bangere  P. Purnaprajna
 \endauthor
 \abstract{In this article we classify quadruple Galois canonical covers
$\varphi$ of singular surfaces of minimal degree. This complements
the work done in \text{\GPsmooth \ns,} so the main output of both
papers is the complete classification of quadruple Galois
 canonical covers of surfaces of minimal degree, both singular and smooth.
 Our results show that the covers $X$ studied in this article
are all regular surfaces and
form a bounded family in terms of geometric genus $p_g$.
 In fact, the geometric genus of $X$ is bounded by $4$. Together with
the results of Horikawa and Konno for double and triple
 covers, a striking numerology emerges that motivates some
 general questions on the existence of higher degree canonical covers.
 In this article, we also answer some of these questions. The
 arguments to prove our results include a delicate analysis
 of the discrepancies of divisors in connection with the ramification 
 and inertia groups of $\varphi$.}
 \endabstract
  \endthanks
 \address{Francisco Javier Gallego: Dpto. de \'Algebra,
  Facultad de Matem\'aticas,
  Universidad Complutense de Madrid, 28040 Madrid,
 Spain}\endaddress
 \email{gallego\@mat.ucm.es}\endemail
 \address{ B.P.Purnaprajna:
 405 Snow Hall,
   Dept. of Mathematics,
   University of Kansas,
   Lawrence, Kansas 66045-2142}\endaddress
 \email{purna\@math.ukans.edu}\endemail
 \thanks{The first author was partially supported by MCT project
number BFM2000-0621. He  is grateful for the hospitality of the
Department of Mathematics of the University of Kansas at
Lawrence. 
The second author is grateful to NSA for supporting this
research project. He is also grateful for the hospitality of the
Departamento de \'Algebra  of the Universidad Complutense de Madrid.}
\endthanks

 \endtopmatter
 \document
 
 \vskip .2 cm
 
 \headline={\ifodd\pageno\rightheadline \else\leftheadline\fi}
 \def\rightheadline{\tenrm\hfil \eightpoint CLASSIFICATION OF QUADRUPLE
   GALOIS CANONICAL COVERS II
  \hfil\folio}
 \def\leftheadline{\tenrm\folio\hfil \eightpoint F.J. GALLEGO \&
   B.P. PURNAPRAJNA \hfil}

 \vskip .3 cm
 
 \heading Introduction \endheading

Canonical covers of surfaces of
minimal degree play a crucial role 
in a variety of contexts including classification problems, the study
of the generation of the canonical 
ring, the study of linear systems on threefolds and the so-called
mapping of the geography of surfaces of general type (see \Pu for a detailed 
motivation.) Among them, Galois canonical covers of degree $4$ are
especially relevant for they behave very differently from both 
canonical double covers and canonical triple covers, as we showed in 
\GPsmooth \ns.  

\smallskip

The classification of double canonical covers of surfaces of minimal
degree was done by Horikawa (see \Hoone \ns). Canonical covers of
degree $3$ were classified by Horikawa (see \text{\Hothree \ns)} 
 and Konno (see \Kon \ns). In \GPsmooth \ns, we
have classified quadruple Galois canonical covers of smooth surfaces 
of minimal degree.
In this article we complete the classification of quadruple Galois
canonical covers of surfaces of minimal degree by studying those
covers whose image  
is a singular surface. We summarize the classification obtained in the
next theorem. Before stating it we need an auxiliary construction: 

\smallskip

Let $X$ be a canonical surface whose canonical bundle is
base-point-free, let $W$ be a singular rational
normal scroll  and let $Y @> q >> W$ be the minimal desingularization
of $W$.
Let $X @> \varphi >> W$ be the canonical morphism
of $X$ (a canonical cover of $W$, for short). 
There exists the following commutative square:

$$ \matrix \overline X & @> \overline q >> & X \cr
  @VV p V &  \hskip -.8 truecm @VV\varphi V \cr
  Y & @> q >>& W \cr
 \endmatrix \ \ (*) $$
 where $\overline X$ is the normalization of
 the reduced part 
of $X \times_W Y$, which is irreducible, and $p$ and $\overline q$ are
 induced by the projections from the
 fiber product onto each factor. Now we can state the following
 result, which inside the article is split up into  
\classthms \ns:

\proclaim{\summary} If $X$ is a canonical surface with base-point-free
canonical bundle, $W$ is  a
singular surface of minimal degree, $X @> \varphi >> W$ is the
canonical 
morphism of $X$ and $\varphi$ is Galois of degree $4$ with Galois
group $G$, then 
$W=S(0,2)$, $X$ is regular  and there exists
a commutative diagram like (*), where $Y @> q >> W$ is the minimal
desingularization of $W$, 
the (normal) surface $\overline X$ has at
worst canonical singularities, the morphism $\overline q$ 
is the morphism from $\overline X$ to its canonical model $X$ and,

\item{1)}
if $\overline q$ is crepant and  $G=\bold Z_2 \times \bold Z_2$, then
$\overline X$ is   
the product over $Y$ of two double covers branched along divisors
$D_2$ and $D_1$ which are linearly equivalent to $2C_0+6f$ and
$4C_0+6f$ respectively and $p$ is the natural morphism from the fiber
product to $Y$; 

\medskip

\item{2)} if $\overline q$ is crepant and $G=\bold Z_4$, 
then $p$ is the composition of two
double covers  
$\overline X_1 @> p_1 >>Y$ branched along a divisor $D_2$ linearly
equivalent to $4C_0+6f$ and $\overline X @> p_2 >> \overline X_1$,
branched along the ramification of $p_1$ and $p_1^*D_1$, with $D_1$
linearly equivalent to $3f$, and with trace zero module $p_1^*\Cal
O_Y(-\frac 1 2 D_1 -\frac 1 4 D_2)$;

\medskip

\item{3)} if $\overline q$ is noncrepant and $G=\bold Z_2 \times
\bold Z_2$,  then $\overline q$ is the blowing-up of $X$ at two smooth
points, $\overline X$ is the
normalization of the fiber 
product over $Y$ of two double covers of $Y$ each branched along a
divisor linearly equivalent to $4C_0+6f$, and $p$ is the natural map
from the normalization of the fiber product to $Y$;  

\medskip

\item{4)} if $\overline q$ is noncrepant and $G=\bold Z_4$, then 
$p$ is the composition of two double covers 
$\overline X_1 @> p_1 >> Y$, branched along a divisor $\Delta_2$, 
and $\overline X @> p_2 >> \overline X_1$,
branched along the ramification of $p_1$ and $p_1^*D_1$ and with trace
zero module $p_1^*\Cal O_Y(-\frac 1 2 (D_1 + C_0) -\frac 1 4 \Delta_2)
\otimes \Cal 
O_{\overline X_1}(\overline C_0)$, where $D_1$ is a divisor on $Y$ and
$\overline C_0 =
p_1^{-1}C_0$,  and either,

\smallskip

\itemitem{}4.1) $D_1 \sim C_0+3f$ and  $D_2 \sim 4C_0+6f$; or 

\smallskip

\itemitem{}4.2)  $D_1 \sim 4C_0+9f$ and $D_2 \sim 2C_0+2f$. 

\medskip  

\noindent Conversely, if $\overline X$ is a normal surface with at worst
canonical singularities, $Y=\bold F_2$ and $\overline X @> p >>
Y$ is either

\smallskip

\item{I.} the fiber product over $Y$ of two double covers branched
  along divisors 
$D_2$ and $D_1$ as described in 1) above, or
the normalization of the fiber product over $Y$ of two
  double covers branched 
  along divisors 
$D_2$ and $D_1$ as described in 3) above; or

\item{II.} the composition of two
double covers  
$\overline X_1 @> p_1 >>Y$, branched along a divisor $D_2$ 
 and $\overline X @> p_2 >> \overline X_1$,
branched along the ramification of $p_1$ and $p_1^*D_1$
 and having trace zero module $p_1^*\Cal
O_Y(-\frac 1 2 D_1 - \frac 1 4 D_2)$, where $D_1$ and $D_2$ are as
described in 2) above; or 

\item{III.} the composition of two double covers 
$\overline X_1 @> p_1 >> Y$, branched along a divisor $\Delta_2$, 
and $\overline X @> p_2 >> \overline X_1$,
branched along the ramification of $p_1$ and $p_1^*D_1$ and with trace
zero module $p_1^*\Cal O_Y(-\frac 1 2 (D_1 + C_0) -\frac 1 4 \Delta_2)
\otimes \Cal 
O_{\overline X_1}(\overline C_0)$, where 
$\overline C_0 =
p_1^{-1}C_0$ and $D_1$ and $\Delta_2$ are as described in 4.1) or
4.2) above, 

\medskip

\noindent
then  
there exists a commutative diagram like (*) 
where $W$ is $S(0,2)$, the morphism $q$ is the minimal
desingularization of $W$, the morphism $\overline q $ is the morphism
from $\overline X$ to its 
canonical model $X$, the morphism 
$\varphi$ is the canonical morphism of $X$ and is
Galois with
Galois group $G$, and $G=\bold Z_2 \times \bold Z_2$ in case I and $G=\bold
Z_4$ in cases II and III. 
\endproclaim

Amidst the landscape of all quadruple Galois canonical covers of
surfaces of minimal degree, covers of singular targets are significant
because among them we find the exceptions to an, otherwise, beautiful
and uniform picture.  The existence of these exceptions adds even more
complexity to the already subtle problem of studying covers of
singular surfaces.

\smallskip

The canonical quadruple Galois covers $X @> \varphi >> W$ of smooth 
surfaces of minimal degree (classified in \GPsmooth \ns) 
and the quadruple Galois covers $\overline X @> p >> Y$ in 1) and 2) of the
above theorem exhibit the same structure.
More precisely, we show that when the canonical covers  
$\varphi$ in \GPsmooth and the covers $p$ of 1) and 2)  
have Galois group $\bold Z_2 \times \bold Z_2$, then the surfaces of general 
type $X$ in \GPsmooth and the surfaces of general type $\overline X$
above 
are always a fiber product of two double covers. Moreover, 
the branch divisors of these double covers also follow a uniform pattern. 
On the other hand, if the Galois group 
is $\bold Z_4$, we show that the morphism $\varphi$ in
\GPsmooth and the morphism $p$ of 1) and 2) are a 
composition of two double covers $p_1$ and $p_2$ such that $p_2$ is
branched along the ramification of $p_1$ and the pullback of a divisor
on the surface of minimal degree. Again the branch divisors of $p_1$
and $p_2$ fit always in the same pattern. 
Thus, after seeing the classification obtained in \GPsmooth and looking
only
at cases 1) and 2) of the present classification, one would be
inclined to  
conjecture this: a surface $X$ (or a closely related birational model
of $X$,  
obtained by a crepant, partial resolution of singularities) 
which is a quadruple Galois canonical cover
of a surface of minimal degree is always either a fiber product of
two double covers or a composition of two double covers $p_1$ and
$p_2$ branched as described above. Cases 3) and 4) of the previous
theorem are exactly the counterexamples to this tempting conjecture.

\smallskip

In cases 3) and 4), the morphism $\overline
q$ is non-crepant, in constrast with cases 1) and 2). In case 3), 
where the Galois group is $\bold Z_2 \times \bold Z_2$, the surface
$\overline X$ is not a fiber product but  the normalization of a fiber
product. The fiber product of the two double covers 
is non-normal precisely because it has a 
double curve that eventually contracts to the vertex $w$ of $W$. 
In case 
4), where the Galois group is $\bold Z_4$, the morphism $p$ is still 
a composition of two double covers, but the cover $p_2$ is branched 
along a divisor which is not a pullback from $Y$. 
The 
main  philosophical reason
why these two exceptions occur is because 
the canonical divisor of $\overline X$ has a fixed part
which contracts eventually to $w$.
  
\smallskip

Another unusual fact that emerges from
the classification of quadruple Galois covers 
is that cyclic quadruple canonical covers of surfaces 
of minimal degree are never simple cyclic. Non-simple cyclic covers are 
not a common phenomenon for surfaces, so its existence in this context is 
interesting.

\smallskip

In this article we also construct families of examples to show the
existence of all the cases that appear in the classification.
We carry out as well a more detailed study of the
singularities of $X$  (see \numbsing \ns). If $w$ is the vertex of
$W$, we see that $\varphi^{-1}\{w\}$ consists of one point 
(cases 1), 2) and 4) of the above theorem) or two points (case 3)).
In cases 3) and 4) the point or points lying over $w$ are smooth, i.e., 
$\varphi$ ``unfolds'' completely the singularity at $w$. 
In the remaining 
cases the singularity over $w$ stays the same or worsens: in case 1) the
point lying over $w$ is an $A_k$ singularity ($A_1$ at best) and in case 
2) it is a $D_4$ singularity. The behavior of the complement 
of the fiber of $w$ is like the behavior of
the canonical covers of smooth surfaces of minimal degree studied in
\GPsmooth \ns: if $G=\bold Z_2 \times \bold Z_2$, one can find covers
for which the complement is smooth, and
if $G=\bold Z_4$, the complement is necessarily singular, having at best $A_1$
singularities. Putting all of the above together we see that 
there do exist smooth quadruple Galois canonical covers of singular surfaces of
minimal degree, but they necessarily belong to case 3). We also
show that all cyclic quadruple canonical covers of surfaces of
minimal degree (smooth or singular) are singular.  

\smallskip

The results in this article show that quadruple Galois
covers of singular surfaces of minimal degree
form a bounded family in terms of both their geometric genus and their
irregularity. 
In fact the 
classification results here show that $p_g\leq 4$ and $q=0$. Together with 
the results of Horikawa and Konno for double and triple covers, 
the following striking numerology emerges for surfaces of general type that 
are Galois canonical covers of singular targets:

$$\matrix 
\text{If deg }\varphi=2, \text{ then } p_g\leq 6, q=0; \cr
\text{if deg }\varphi=3, \text{ then } p_g\leq 5, q=0;  & \ns \ns
\text{and} \cr 
\text{if deg }\varphi=4, \text{ then } p_g\leq 4, q=0 & \ns \ns \ns
\ns \ns \ns \ns .\cr 
\endmatrix$$

Since the smallest projective space containing a singular scroll is
$\bold P^3$, this pattern suggests that there 
do not exist higher degree canonical covers of
singular rational normal scrolls, so we pose the following

\smallskip

\noindent{\bf Question:} Let $X @> \varphi >> W$ be a canonical cover 
of a singular surface of minimal degree $W$. Is deg$\varphi \leq 4$?

\smallskip

There are strong hints towards a positive solution to the question above: 
in \GPsmooth \ns, Corollary 3.3, we prove that there are no regular Galois 
canonical covers of degree prime $p \geq 5$ of a surface of minimal 
degree $W$, smooth or singular. The significance of our question 
becomes clear once we realize the following: if the answer is positive, 
then, having in account our previous results for odd degree covers 
(see \GPtrans \ns, Corollary 3.2), there will be no canonical 
covers of degree odd bigger than $3$ of surfaces of minimal degree, 
except perhaps covers of $\bold P^2$.

\smallskip

As it often happens in classification problems, the special cases are
not  easier to deal with. 
A posteriori we see that quadruple
Galois canonical covers of singular surfaces of minimal degree
represent a smaller portion if we compare them with the covers of
smooth surfaces. However, the difficulties to carry out the
classification in the singular target case are much greater. 
We glimpse
them by briefly commenting on the strategy we follow and the
techniques we employ. To study the cover $X @ > \varphi >> W$ of a
singular scroll $W$ the first thing we do is to ``desingularize''
$\varphi$ using the commutative diagram (*). Once this is done, the
only information avalaible on $\overline X$ is that, by construction,
$\overline X$ is a normal, irreducible surface. Likewise, little is
known of $\overline q$, just that it is a birational morphism between
$\overline X$ and $X$. At this point, the best possible situation one
can hope in order to continue the study of $X$ is that $\overline X$
have canonical singularities and $\overline q$ be crepant, for in such
case one can expect to deal with $\overline X @ > p >> Y$ in much the
same way as with a canonical cover $X @> \varphi >> W$ of a smooth
surface $W$. Thus, the crux of our argument, which is contained in the
proofs of  \singlemmas \ns, is to find out if this
favorable situation happens always, or, if not, if we can control the
``badness'' of $\overline X$ and $\overline q$. To settle the question
we have to study the possible discrepancies of $\overline q$ and the
inertia groups of the ramification of $p$. It finally turns out, as the
reader already knows, that $\overline q$ is not always crepant, but in
the case it is not, by the work done in  \singlemmas  \ns, we are able to
narrow the field and say that 
$\overline q$ has to fulfill very concrete specifications (see
\text{\singlemma \ns,} 2) and
\text{\thmsinglept \ns,} 2)). This, after still some more involved work,
especially when $\overline q$ is not crepant,
makes the problem tractable at the end. Likewise, after \singlemmas
 we find that $\overline X$ is not too bad either (it 
has at worst rational, $2$-Gorenstein singularities). This part of
the tale has an even happier ending, since we eventually prove that the
singularities of $\overline X$ are indeed canonical. 

 \medskip

{\bf Acknowledgements:} We are very thankful to Dale Cutkosky and
N. Mohan Kumar for patiently listening to some of our arguments and 
for their insightful and helpful comments. We 
are also very grateful to Rita Pardini for her valuable comment on 
non-simple cyclic covers. Finally, we thank the referee of the 
announcement \CR of these results, whose advice helped us to improve 
the exposition of the present article.

 \heading 1. Background material
 \endheading

\noindent {\bf Convention:} 
We work over an algebraically closed field of
characteristic $0$. 

\medskip
 
\noindent{\bf{\notation}}\ns : 
{\it We will follow these conventions:}

\smallskip
 
\noindent  1) Throughout this article, unless otherwise stated,
$W$ will be an embedded
projective algebraic surface which is a
cone over a smooth rational normal curve. Thus $W$ has 
minimal degree, for its degree is equal to its codimension in
projective space plus $1$. If the rational normal curve has degree $e$
we will denote $W$ as $S(0,e)$. We will denote by $w$ the vertex of
$W$. 
 
 \smallskip
 
\noindent  2) Throughout
this article, unless otherwise stated, 
$X$ will be a projective algebraic normal
 surface 
 with at worst canonical singularities. 
 We will denote by  $\omega_X$ the canonical bundle of $X$.

\smallskip
 
 \noindent We recall the following standard notation: 
 
\smallskip

\noindent  3) Let $e \geq 2$. 
By $\bold F_e$ we denote the Hirzebruch surface whose
minimal section have self intersection $-e$. Let $C_0$
denote the minimal section of $\bold F_e$ and let $f$ be one of the
fibers of $\bold F_e$. Recall that $S(0,e)$
is the
image of $\bold F_e$ by the complete linear series
$|C_0+ef|$.

\proclaim{\defcan} Let $X$ and $W$ be as in the previous notation. We
will say that a surjective morphism $X @> \varphi >> W$ is a canonical
cover of $W$ if $X$ is surface of general type whose canonical
bundle $\omega_X$ is ample  
and
base-point-free and
$\varphi$ is the canonical morphism of $X$. 
\endproclaim

In this paper we study Galois covers $\varphi$ of $W$. Since $W$ is
singular, $\varphi$ is not in general flat. 
However the strategy will be to study an auxiliary, flat Galois
cover. We recall here  some known or easy facts regarding the algebra
structure associated to a flat,
quadruple Galois cover. For proofs of \Galquad see  \Ca \ns, \HM \ns and \Pa
\ns, and also \GPsmooth \ns, Proposition 2.4; for \biproconv and
\Galcycliconv \ns, see  
\GPsmooth \ns, Propositions 2.9 and 2.10.

\proclaim{\Galquad } 
Let $\frak X$ and $\frak Y$ be two algebraic 
varieties and let $\frak X @> p >> \frak Y$ be a flat,
Galois cover of degree $4$.

 \item{1)} If $G=\bold Z_4$, then  $p_*\Cal O_\frak X$ splits as
 $$p_*\Cal O_\frak X = \Cal O_\frak Y \oplus L_i^* \oplus L_{-1}^* \oplus
 L_{-i}^*$$ where the line bundles on $\frak Y$,   
 $\Cal O_\frak Y, L_i^*, L_{-1}^*$ and $ L_{-i}^*$ are the eigenspaces of
 $1,i,-1$ and $-i$ respectively.

 There exist effective Cartier divisors $D_{11}, D_{12}, D_{23}$ and $D_{33}$
 on $\frak Y$ such that 
 $D_{11}+D_{23}=D_{12}+D_{33}$ and the following 
 $$\matrix L_i \otimes L_i &=& L_{-1} \otimes \Cal O_{\frak Y}(D_{11}) \cr
 L_i \otimes L_{-1} & = & L_{-i} \otimes \Cal O_{\frak Y}(D_{12})\cr
 L_i \otimes L_{-i} & = & \Cal O_{\frak Y}(D_{11}+D_{23})  
 \cr
 L_{-1}\otimes L_{-1} & = & \Cal O_{\frak Y}(D_{12}+D_{23})\cr
 L_{-1}\otimes L_{-i} & = & L_i \otimes \Cal O_{\frak Y}(D_{23})\cr
 L_{-i} \otimes L_{-i} & = & L_{-1} \otimes \Cal O_{\frak Y}(D_{33})\cr 
 \endmatrix 
 \qquad \Galquadone    
 $$
 and the multiplicative structure of  $p_*\Cal O_\frak X$ is as
 follows: 
 
 $$\displaylines{L_i^* \otimes L_i^* @> \cdot D_{11} >> L_{-1}^* \cr
 L_i^* \otimes L_{-1}^* @> \cdot D_{12} >> L_{-i}^*\cr
 L_i^* \otimes L_{-i}^* @> \cdot D_{11}+D_{23}  >> \Cal O_\frak Y  \cr
 L_{-1}^* \otimes L_{-1}^* @> \cdot D_{12}+D_{23} >> \Cal O_\frak Y \cr
 L_{-1}^*\otimes L_{-i}^* @> \cdot D_{23} >> L_i^*\cr
 L_{-i}^* \otimes L_{-i}^*  @> \cdot D_{33} >> L_{-1}^* \ . \cr}$$

 \smallskip
 
 \item{2)} If $G=\bold Z_2 \times \bold Z_2$, 
 then $p_*\Cal O_\frak X$ splits as
 $$p_*\Cal O_\frak X = \Cal O_\frak Y \oplus L_1^* \oplus L_2^* \oplus
 L_3^* \ ,$$ where $\Cal O_\frak Y$, $L_1^*$,
   $L_2^*$ and $L_3^*$ are eigenspaces
 and there exist effective Cartier divisors 
   $D_1$, $D_2$ and $D_3$ such that 
 $L_i^{\otimes 2}=\Cal O_{\frak Y}(D_j+D_k)$ and 
 $L_j \otimes L_k = L_i \otimes \Cal O_{\frak Y}(D_i)$ with 
$i \neq j$, $j \neq k$ and $k 
 \neq i$, and the multiplicative structure of $p_*\Cal O_\frak X$ is as
 follows: 
 $$\displaylines{L_i^* \otimes L_i^* @> \cdot D_j+D_k >> \Cal O_\frak Y \ ,
   \cr
 L_j^* \otimes L_k^* @> \cdot D_i >> L_i^*
 \ .}$$ 
\endproclaim

 \proclaim{\biproconv} 
 Let $\frak X$ and $\frak Y$ be normal algebraic varieties.
 \item{1)} If 
 $\frak X @> p >> \frak Y$ 
 is the natural map onto $\frak Y$ 
from the fiber product
 over 
 $\frak Y$ 
 of two flat double covers $\frak X_1 @> p_1 >> \frak Y$ and $\frak X_2 @>
 p_2 >> \frak Y$, then $p$ is a Galois cover with Galois group $\bold
 Z_2 \times \bold Z_2$. 
 \item{2)} If in addition  
 $L_2^*$ and $L_1^*$ are the trace zero 
 modules of $p_2$ and $p_1$ respectively, then 
 $$p_*\Cal O_\frak X=\Cal O_\frak Y \oplus L_1^* \oplus L_2^* \oplus
 (L_1^* \otimes L_2^*) $$ and, if $\frak Y$ is locally Gorenstein,
 then  $\frak X$ is locally Gorenstein
 and $\omega_\frak X= p^*(\omega_\frak Y \otimes L_1 \otimes L_2)$. 
 \endproclaim

 \proclaim{\Galcycliconv} 
 Let $\frak X$ and $\frak Y$ be normal algebraic varieties. 
 \item{1)} If 
 $\frak X @> p >> \frak Y$ 
 is the composition of a flat double cover $\frak X' @> p_1 >>
 \frak Y$ branched along a divisor $D_2$,  
 followed by a flat double cover $\frak X
 @> p_2 >> \frak X'$, 
 branched along  the ramification
 locus of $p_1$ and $p_1^*D_{1}$, where $D_1$ is a divisor on $\frak
 Y$,  
 then $p$ is a Galois cover with Galois group $\bold
 Z_4$. 
 \item{2)} If in addition $L_2^*$ is the trace zero module of $p_1$  and 
 $p_1^*L_1^*$ is the trace zero module of $p_2$, then 
 $$p_*\Cal O_\frak X=\Cal O_\frak Y \oplus L_1^* \oplus L_2^* \oplus
 (L_1^* \otimes L_2^*) $$ and, if 
 $\frak Y$ is locally Gorenstein, then $\frak X$ is
 locally Gorenstein and 
 $\omega_\frak
 X=p^*(\omega_\frak Y
 \otimes L_1 \otimes L_2 ).$
 \endproclaim

\heading 2. The desingularization diagram \endheading

The covers we want to describe and classify in this article are
Galois canonical covers $\varphi$ of a singular scroll $W$. These covers are
finite but, precisely because $W$ is singular, they are not in general
flat. Flat covers are more tractable though, since their associated
algebra structure is locally free. Thus the first thing we do is to
``desingularize'' $\varphi$, that is, we will ``make'' $W$ smooth and
$\varphi$ flat. We
construct the following {\it desingularization diagram}  for $\varphi$:

\proclaim{\diagram} Let $X @> \varphi >> W$ be a canonical
cover and let $Y @> q >> W$ be the
minimal desingularization of $W$. We define $\overline X$ as the
normalization of 
the reduced part 
of 
$X \times_W
Y$. The surface $\overline X$ is irreducible 
and 
fits in the following commutative
diagram:
 $$\matrix \overline X & @> \overline q >> & X \cr
         @VV p V &  \hskip -.8 truecm @VV\varphi V \cr
 Y & @> q >>& W \cr 
 \endmatrix \  \diagramone$$ 
 where   
 $p$ and $\overline q$ are
 induced by the projections from the fiber product onto each factor.
\endproclaim

For our purposes, the key point of the above construction is that if one of 
$\varphi$ or $p$ is Galois with given Galois group $G$, so is the
other: 

\proclaim{\diagramGal } 
Let
$X @> \varphi >> W$ be a canonical cover and let $\overline X @> p >> Y$
be as in \diagramone \ns. Then $\varphi$ 
is a Galois cover with Galois group $G$ if and only if 
$\overline X @> p >> Y$ is a 
Galois cover with Galois group $G$. 
\endproclaim

{\it Proof.} We first assume that $\varphi$ is Galois and its Galois
group is $G$. 
Since $\Cal O_X$ is integral over $\Cal O_W$, 
$\Cal O_{\overline X}$ is also an integral extension of $\Cal O_Y$.  
 By construction $\overline X$ is normal. 
 Recall that $W$ is a cone over a smooth rational normal curve and let 
$w$
 be the singular point of $W$. Then $\overline X$ and $X$ are
 birational, in fact, isomorphic outside the points of $X$ lying over
 $w$. Therefore 
 $\Cal O_{\overline X}$ is in fact the integral closure of $\Cal
 O_Y$ in $\Cal K(X)$, so  
 $\overline X @> p >> Y$ is
also
 a Galois cover with the same Galois group $G$. 
The argument to show the converse is analogous. \qed

\bigskip

 We state now two lemmas about $p$ and $\overline q$. The first of
 them recalls well-known facts on rational singularities, so we state
 it without a proof:

 \proclaim{\ratlemma} Let $\overline{\frak X}@> \overline q >>
 \frak X$ be a
 birational morphism between two normal surfaces. 
 If $\frak X$ has rational singularities, then
 \item{1)} $\overline{\frak X}$ also has rational singularities, and
 \item{2)} every reduced cycle of $\overline{\frak  X}$ contracted to a point
   by $\overline q$ has arithmetic genus $0$.  
 \endproclaim

 \proclaim{\discrlemma} With the notation of \diagram \ns, if 
 $X @> \varphi >> W$ 
is a Galois cover with group $G$, then

 \item{1)} On $X$ and $\overline X$ there exist canonical divisors
$K_X$ and $K_{\overline
   X}$ which are $G$-invariant and such that 
 $$K_{\overline X} \equiv {\overline q}^*K_X + a(F_1 + \cdots + F_k),
  \quad \prsinglemmaone 
$$
 where $\equiv$ means numerical equivalence, $a$ is a nonnegative
 rational number and $F_1, \dots, F_k$ are the components of the
 exceptional locus of $\overline q$.
 
 \item{2)} If in addition $\overline X$ is locally Gorenstein, then there
   exist $K_X$ and $K_{\overline
   X}$ as above and such that 
 $$K_{\overline X} = {\overline q}^*K_X + a(F_1 + \cdots + F_k),
 \quad \prsinglemmatwo 
$$
 where $a$ is a nonnegative integer.  
 
 \endproclaim
 
 {\it Proof.} 
 We consider the exceptional locus of $\overline q$. Recall that
 $Y$ is a Hirzebruch surface and let $C_0$ be its minimal section. 
 Any curve
 $F_i$ in
 the exceptional locus of $\overline q$  maps onto $C_0$ by $p$. 
 
 Since the cover $p$ is Galois by \diagramGal \ns, 
$G$ acts transitively on the set 
 $\{F_1, \dots, F_k\}$. Recall also that $X$ and
 $\overline X$ are both normal, $X$ has at worst canonical
 singularities (in particular $X$ is also locally Gorenstein) and, by \ratlemma
 \ns,
 $\overline X$ has at worst rational singularities (in particular,
 $\overline X$ is locally $\bold Q$-Gorenstein). 
Then 
 one can find $G$-equivariant canonical divisors $K_X$ and
 $K_{\overline X}$ such that 
 $$K_{\overline X} \equiv {\overline q}^*K_X + a(F_1 + \cdots + F_k) $$
 Then $a$ is a nonnegative rational number, because $X$ has
 canonical singularities. This proves 1)
 
 \smallskip
 
 If, in addition, $\overline X$ is locally Gorenstein
 in the previous formula we can write equality instead of numerical
 equivalence and 
 $a$ is an integer, for both $K_{\overline X}$ and $K_X$ are
 Cartier divisors. This proves 2). \qed

\bigskip

The philosophy we follow now is this: instead of describing directly
the quadruple Galois canonical covers $X @> \varphi >> W$ of $W$ 
we will describe Galois covers
$\overline X @> p >> Y$. We will describe also in a precise way the
relation between
$\overline X$ and $X$, that is, how one passes from $\overline X$
to $X$ and viceversa. That means to describe
the morphism $\overline q$.  We split the study of $\overline q$ in
two theorems, 
 depending on whether $\{\varphi^{-1}(w)\}$ consists of one or several
 points.

 \proclaim{\singlemma } Let $W$ be a singular rational 
 normal scroll and let $X @> \varphi >> W$ be a canonical cover. 
 Let $w$ be the singular point of $W$ and let $\overline X, Y, p, q$ and
 $\overline q$ be as in \diagram \ns.
  If $X @> \varphi >> W$ is Galois of degree $4$ 
 and $\{\varphi^{-1}(w)\}$ is not a single point,
 then $\overline X$ 
has at
 worst canonical 
 singularities and one of the following happens:
 
 \item{1)} either $\overline q$ is crepant (i.e., $K_{\overline
   X}=\overline q^*K_X$); or
 
 \item{2)}  $W=S(0,2)$, $\varphi^{-1}\{w\}$ consists of two smooth points $x_1$
   and $x_2$ and  $\overline{X} @> {\overline q} >> X$ is the blowing
   up of $X$ at $x_1$ and $x_2$. 
 \endproclaim

 {\it Proof.}
 Recall that $Y$ is a Hirzebruch surface $\bold F_e$ with $e
 \geq 2$, that $C_0$ is its minimal section and that we call $F_1,
 \dots F_k$ the irreducible components of the exceptional locus of
 $\overline q$, which are mapped each onto $C_0$ by $p$. Let $G$ be the
 Galois group of $\varphi$.  
 Since the order of $G$ is $4$ and the cardinality of $\varphi^{-1}\{w\}$
 is greater than one, 
 the cardinality of $\varphi^{-1}\{w\}$ is in fact $2$ or $4$. We treat
 these two cases separately: 
 
 \smallskip
 
 \noindent {\it Case 1:} Cardinality of $\varphi^{-1}\{w\}$ equals $4$. 
In this case $\varphi$ is \'etale at $w$, and hence is flat, 
so $X \times_W Y$ is the blowing up of $X$ at the four points  $x_1,
 \dots, x_4$, lying over $w$. 
Thus $X \times_W Y$ is irreducible, reduced and normal, so, in this
case, $\overline X = X \times_W Y$.
 Since $\varphi$ is \'etale at an analytic neighborhood of $w$, the
 singularities at $x_1,\dots, x_4$ are analytically isomorphic to the
 singularity at $w$, i.e., they are all $A_1$ singularities. Thus
 $\overline q$ resolves $x_1, \dots, x_4$. Therefore 
$\overline q$ is crepant,
 i.e., $K_{\overline X}=\overline q^*K_X$. 
 Finally, since  
$X$ has at worst canonical singularities, so does $\overline X$ (the
canonical singularities of $\overline X$
correspond  to the singular points
of $X$ different from $x_1, \dots, x_4$).
 \smallskip
 
 \noindent{\it Case 2:} Cardinality of $\varphi^{-1}\{w\}$ equals
 $2$. We call 
 $x_1$ and $x_2$ the points in $\varphi^{-1}\{w\}$. 
 Given a subgroup $G'$ of $G$,
 let $X'$ be the quotient of $X$ by $G'$. We
 also have a way of decomposing $\varphi$, namely, 
 $$X @> \varphi_2 >> X' @> \varphi_1 >> W \ ,$$ 
 where $X'$ is normal and $\varphi_1$ and $\varphi_2$ are Galois
 covers with Galois group $G/G'$ and $G'$ respectively.
  
 Let $\overline{X'}$ be the normalization of the reduced part 
of the fiber product of $X'$ and $Y$
 over $W$ and let $\overline{X'} @> \overline q >> X'$ and
 $\overline{X'} @> p_1 >> Y$ be the projections to each factor of the
 product. 
As it happened with $\overline X$, $\Cal O_{\overline{X'}}$ is the
 integral closure of $\Cal O_Y$ inside $\Cal K(X')$ and therefore $p_1$
 is a Galois cover with group $G/G'$. Now let $\overline{X''}$ be so
 that $\Cal O_{\overline{X''}}$ is the integral closure of $\Cal
 O_{\overline{X'}}$ in $\Cal K(X)$ and let $p_2$ be the morphism
 induced between $\overline{X''}$ and $\overline{X'}$.  Then $\Cal
 O_{\overline{X''}}$ is 
 normal and
 integral over $\Cal O_Y$, hence it is the integral closure of $\Cal
 O_Y$ in $\Cal K(X)$. But by construction, so is $\Cal
 O_{\overline{X}}$, hence $\overline{X}=\overline{X''}$ and 
 we get the following commutative diagram: 
 
 $$ \matrix \overline X & @> \overline q >> & X \cr
 @VV p_2 V &  \hskip -.8 truecm @VV\varphi_2 V \cr
 \overline {X'} & @> q' >> & X' \cr
         @VV p_1 V &  \hskip -.8 truecm @VV\varphi_1 V \cr
 Y & @> q >>& W \cr 
 \endmatrix \  \diagramtwo $$  
 with $p=p_1 \circ p_2$ and $\varphi=\varphi_1 \circ \varphi_2$. 
 Let now $G'=\text{Stab}x_1=\text{Stab}x_2$.
 Let $x_1'=\varphi_2(x_1)$ and 
  $x_2'=\varphi_2(x_2)$. Since $G'$ is the stabilizer of both $x_1$ and
  $x_2$, then $x_1' \neq x_2'$ and $\varphi_1$ is \'etale at a
  neighborhood of $w$. By the same
  arguments as in \text{Case 1} $\overline{X'} @> q' >> X'$ 
  is the blowing up of $X'$ at
  $x_1'$ and $x_2'$. 
 Then $p_1$ is also \'etale on an analytic neighborhood of $C_0$. Let
  $E_1$ and $E_2$ be the exceptional divisors of $q'$. Then
  $p_1^*C_0=E_1+E_2$, $E_i$ is isomorphic to $\bold P^1$, 
 $\overline{X'}$ is smooth at every point of an
  analytic neighboorhood $U$ of $E_1$
  and $E_2$, $E_1 \cdot E_2 =0$ and $E_1^2=E_2^2=C_0^2=-e$. Now we examine the
  singularities of $\overline X$.  
 By construction the exceptional locus of ${\overline q}$ 
 is mapped by $\overline q$ on $x_1$
  and  $x_2$ and is in fact $p^{-1}(C_0)$. Thus $\overline X -
  p^{-1}(C_0)$ and $X - \{x_1,x_2\}$ are isomorphic, so $\overline X$
  has at worst canonical singularities outside $p^{-1}(C_0)$. 
 By construction $\overline X$ is normal, hence locally Cohen-Macaulay, and
  since $U$ is smooth, $p_2$ is a flat, degree $2$ morphism 
  when restricted to $p_2^{-1}(U)$. Therefore $p_2^{-1}(U)$ is 
  locally Gorenstein. 
 By \ratlemma  $\overline X$ has rational singularities, hence
 $p_2^{-1}(U)$ has Gorenstein rational singularities, i.e., canonical
 singularities. This proves that $\overline X$ has canonical
 singularities. 
 
 \smallskip
 
 Now we study $\overline q$. For that  
 we study how $p_2$ is at
 $E_1$ and $E_2$. Recall that $p_2$ is a double cover branched along a
 divisor of $\overline{X'}$. Since $G$ acts transitively, there are
 only two possibilities for $E_1$ and $E_2$: either $E_1$
 and $E_2$ are both in the branch locus of $p_2$ or none of them are. 
 Now we deal with the two possibilities. First, let us assume that
 neither $E_1$ nor $E_2$ are in the branch locus of $p_2$. Let
 $F_i=p_2^*E_i$. Then $F_i$ is a Cartier divisor in $\overline X$ and it
 is a reduced curve, $F_1 \cdot F_2 =0$ and $F_i^2=-2e$. Since $X$ has
 rational singularities $F_i$ has arithmetic genus $0$ by \ratlemma \ns.
 
 Recall that we have shown $\overline X$ is locally Gorenstein. Then from 
 adjunction and from \prsinglemmatwo  we get
 $$-2 = (K_{\overline X} + F_i)\cdot F_i= ({\overline q}^*K_X + aF_1 +
 aF_2 + F_i) \cdot F_i = -2e(a+1) \ ,$$
 with $a$ a nonnegative integer. 
This leads to a contradiction, for 
 $e$ is an integer greater than or equal to $2$.  
 
 Then the only possibility left is that both $E_1$ and
 $E_2$ are in the branch locus of $p_2$. Now let
 $F_i=p_2^{-1}E_i$. Then $F_i$ is isomorphic to $E_i$ and therefore to
 $\bold P^1$ and $2(F_1+F_2)=p^*_2(E_1+E_2)=p^*C_0$. Using
 \prsinglemmatwo and the commutativity of diagram \diagramtwo we obtain
 $$\displaylines{K_{\overline X}={\overline q}^*K_X + a(F_1 + F_2)= \cr
   {\overline 
   q}^*\varphi^*K_X + a(F_1 + F_2)=p^*(C_0+ef)+a(F_1+F_2)  \ \ 
 \prsinglemmathree }$$
 with $a$ a nonnegative integer. 
 On the other hand, let us denote by $R$ the ramification divisor of
 $p$. Then we have the formula
 $$K_{\overline X}=p^*K_X + R \sim p^*(-2C_0-(e+2)f) + R  \ \
 \prsinglemmafour $$
 where $\sim$ means linear equivalence. Thus from \prsinglemmathree and
 \prsinglemmafour  
 we obtain
 $$R \sim p^*(3C_0+(2e+2)f)+a(F_1+F_2) \ .  \ \ \prsinglemmafive $$
 Since $E_1$ and $E_2$ are in the branch locus of $p_2$, $F_1$ and
 $F_2$ are in the support of $R$. Since $\overline X$ is normal, the
 multiplicity of $E_1$
 and $E_2$ in the branch locus of $p_2$ is $1$. Recall also that $p_1$
 is \'etale at $E_1$ and $E_2$. Thus the multiplicity of $F_1$ and
 $F_2$ in $R$ is also $1$ and we can write $R = R_1 + F_1 + F_2$, where
 $R_1$ is an effective divisor not containing $F_1$ or $F_2$. Thus 
 $R_1 \cdot F_i \geq 0$ and by \prsinglemmafive 
 $$R_1 \sim p^*(2C_0+(2e+2)f) + (a+1)(F_1 + F_2) \ .$$ Putting these two
 pieces of information together we get
 $$\displaylines{  0 \leq R_1 \cdot (F_1 + F_2) = \frac 1 2 p^*(2C_0 +
   (2e+2)f) \cdot 
 p^*C_0 + \frac 1 4 (a+1)(p^*C_0)^2 = \cr 
 2(2C_0 +
   (2e+2)f) \cdot C_0 + (a+1)C_0^2 = 4 -e(a+1) \ \ \prsinglemmasix }$$
 
 Recall that $e$ is an integer greater than or equal to $2$. Then from
 \prsinglemmasix we obtain that $2 \leq e \leq 4$ and $a=0,1$. Moreover, 
 if $e=3,4$, then $a=0$ so in this case $K_{\overline X}={\overline q}^*K_X$. 
 If $e=2$, then either $K_{\overline X}={\overline q}^*K_X$ or 
 $K_{\overline X}={\overline q}^*K_X + F_1 + F_2$. In the latter case,
 $$R_1 \cdot (F_1 + F_2) =0 \ \ \ \prsinglemmasixb $$ 
 hence $\overline X$ is smooth at every point of
 $F_1$ and $F_2$, and $F_1^2=F_2^2=-1$. Then by Castelnuovo's
 contractibility criterion $X$ is smooth at
 $x_1$ and $x_2$ and $\overline q$ is the blowing-up of $X$ at $x_1$ 
 and $x_2$. \qed
 
 \bigskip
 
 Now we study the desingularization diagram \diagramone
  when the inverse image of $w$ in $X$ is a single point:

 \proclaim{\thmsinglept}  Let $W$ be a singular rational 
 normal scroll and let $X @> \varphi >> W$ be a canonical cover. 
 Let $w$ the singular point of $W$ 
and let $\overline X, Y, p, q$ and
 $\overline q$ be as in \diagram \ns.
 If $X @> \varphi >> W$ is Galois of degree $4$ 
 and Galois group $G$, and $\{\varphi^{-1}(w)\}$ is a single point,
 then one of the following happens:
 
 \item{1)} either the surface $\overline X$ has at
 worst canonical 
 singularities, $\overline q$ is crepant (i.e., $K_{\overline
   X}=\overline q^*K_X$) and $W=S(0,2)$; or

\item{2)} the surface $\overline X$ has at worst canonical
  singularities, $\overline q$ is crepant and $G=\bold Z_4$; or
 
 \item{3)}  the surface $\overline X$ is locally $2$-Gorenstein, has at
worst rational singularities and $W=S(0,2)$;
 $G$ is $\bold Z_4$ and the exceptional
 divisor of $\overline q$ is a smooth line $F$
with inertia group $G$, $F^2=-\frac 1 2$ and 
$2K_{\overline X} =
 {\overline q}^*2K_X + 4F$
 for suitable canonical divisors $K_{\overline X}$ and $K_X$. 
 
 \endproclaim

 {\it Proof.}
 Let $\varphi^{-1}\{w\}=\{x\}$. We treat separately the cases of
 $G=\bold Z_4$ and $G=\bold Z_2 \times \bold Z_2$. 
 
 {\it Case 1:} $G=\bold Z_4$. 
  Recall the multiplicative structure of
  $p_*\Cal O_{\overline X}$:
   $$p_*\Cal O_{\overline X}= \Cal O_Y \oplus L_{i}^* \oplus L_{-1}^* \oplus
  L_{-i}^* \ ,$$
  where $L_i^*$, $L_{-1}^*$ and $L_{-i}^*$ are eigenspaces for $i$,
  $-1$ and $-i$ 
  respectively and the multiplication is given by divisors 
 $D_{11}, D_{12}, D_{23}$
  and $D_{33}$, as described in \Galquad \ns.
 
\smallskip

We now describe further the branch locus of $p$. 
Recall that $L_{-1}^{\otimes 2}=\Cal O_Y(D_{12}+D_{23})$, and, since
the
  double cover of $Y$ corresponding to $\Cal O_Y
\oplus L_{-1}$ can be
  taken to be normal, we may assume that $D_{12}+D_{23}$ has no multiple
  components.
    On the other hand, $D_{11}+D_{23}=D_{12}+D_{33}$, hence $D_{23}\subset
  D_{33}$ and $D_{12} \subset D_{11}$, so we can write
  $D_{33}=D_{23}+D_{33}'$ and $D_{11}=D_{12}+D_{11}'$, with
  $D_{11}'=D_{33}'$. Thus the multiplicative structure of 
 $p_*\Cal O_{\overline X}$ and
  the relation between the eigenspaces is summarized as follows:
 
  $$\displaylines{L_i^* \otimes L_i^* @> \cdot D_{11}'+D_{12} >> L_{-1}^* \cr
   L_i^* \otimes L_{-1}^* @> \cdot D_{12} >> L_{-i}^*\cr
   L_i^* \otimes L_{-i}^* @> \cdot D_{11}'+D_{12}+D_{23}  >> \Cal O_Y\cr
   L_{-1}^* \otimes L_{-1}^* @> \cdot D_{12}+D_{23} >> \Cal O_Y \cr
   L_{-1}^*\otimes L_{-i}^* @> \cdot D_{23} >> L_i^*\cr
   L_{-i}^* \otimes L_{-i}^*  @> \cdot D_{11}'+D_{23} >> L_{-1}^*
   \cr}$$
\noindent and 
 $$\displaylines{
L_i \otimes L_i = L_{-1} \otimes \Cal O_Y(D_{11}'+D_{12}) \cr
 L_i \otimes L_{-1} = L_{-i} \otimes \Cal O_Y(D_{12})\cr
 L_i \otimes L_{-i} = \Cal O_Y(D_{11}'+D_{12}+D_{23})\cr
 L_{-1}\otimes L_{-1} = \Cal O_Y(D_{12}+D_{23})\cr
 L_{-1}\otimes L_{-i} = L_i \otimes \Cal O_Y(D_{23})\cr
 L_{-i} \otimes L_{-i} = L_{-1} \otimes \Cal O_Y(D_{11}'+D_{23}) \ .\cr}$$
 
  Then the ramification of $p$ falls only onto the components of
$D_{11}'$, 
 $D_{12}$ and
  $D_{23}$. Since $\overline X$ is normal, by computing locally the
  ramification  lying
  over the generic points of each component of $D_{11}'$, $D_{12}$ and
  $D_{23}$, we can conclude that $D_{11}'+D_{12}$ and 
$D_{11}'+D_{23}$ have
no multiple components either.
  
  Now we discuss how $D_{12}$ and $D_{23}$ are. 
We will show that either $D_{12}$ or $D_{23}$ is a multiple of
$C_0$. We will see it by looking at the intersection of $D_{12}$ and
$D_{23}$ with $C_0+ef$. Let $D$ be a smooth irreducible
curve in $|C_0+ef|$ so that the  
pullback of $D$ by $p$ is also a smooth curve $C$ in $|\overline q^*
\varphi^*\Cal O_W(1)|$. The curve $D$ is isomorphic to $\bold
P^1$. Since $D$ does not meet $C_0$, it corresponds 
to a smooth hyperplane section of $W$ which avoids $w$. Thus, by
adjunction, $\omega_C=p^*(q^*\Cal O_W(2) \otimes \Cal O_D)$. Then,
using relative duality and arguing in a similar fashion as in the
proof of \GPsmooth \ns, Proposition 1.3,  
 we conclude that, for some permutation $\tau$ of
$\{i,-1,-i\}$,  $(L_{\tau(i)} \otimes L_{\tau(-1)})|_D=q^*\Cal O_W \otimes
\omega_D^*=L_{\tau(-i)}|_D$.  Thus either $(L_i \otimes L_{-1})|_D =
L_{-i}|_D$ or 
$(L_{-1} \otimes L_{-i})|_D = L_i|_D$ or $(L_i \otimes L_{-i})|_D =
L_{-1}|_D$.  In the
first case we have $D_{12} \cdot D=0$. This implies that $D_{12}$ is a
multiple of $C_0$. Likewise, in the second case we have that $D_{23}$
is a multiple of $C_0$. Finally, if $(L_i \otimes L_{-i})|_D =
L_{-1}|_D$, we see that $2(D_{11}'+D_{12}+D_{23})$ and
$(D_{12}+D_{23})$ have the same restriction to $D$ and hence
$(2D_{11}'+D_{12}+D_{23}) \cdot D=0$. This implies that $D_{11}',
D_{12}$ and $D_{23}$ are all multiple of $C_0$.

\smallskip
 
Therefore we may rename
 $\{L_i, L_{-1}, L_{-i}\}$ as $\{L_1,L_2,L_3\}$ and
 $\{D_{11}',D_{12},D_{23}\}$ as $\{D_1,D_2,D_3\}$ so that
 
 $$\displaylines{L_1 \otimes L_1= L_2 \otimes \Cal O_Y(D_1+D_2) \cr
 L_1 \otimes L_2 = L_3 \otimes \Cal O_Y(D_{2})\cr
 L_1 \otimes L_3 = \Cal O_Y(D_1+D_{2}+D_{3})\cr
 L_2\otimes L_2 = \Cal O_Y(D_{2}+D_{3})\cr
 L_2\otimes L_{3} = L_1 \otimes \Cal O_Y(D_{3})\cr
 L_{3} \otimes L_{3} = L_{2} \otimes \Cal O_Y(D_{1}+D_{3})  \  \cr}$$
 and $D_2$ is a multiple of $C_0$. Thus $D_2$ is either
$0$  or
 $C_0$, since we know that $D_{2}$ has no multiple components. 
 
\smallskip 

We study now the ramification of $p$ and the canonical divisor of
 $\overline X$.  
 The ramification $R_1$ lying over $D_{2}+D_{3}$ has inertia group $\bold
 Z_4$, i.e., the points of $R_1$ have stabilizer $\bold Z_4$.
 To compute the rest of the ramification, we work on
 $U=Y-\{D_{2}+D_{3}\}$, and there it is clear that the only
 ramification lies over $D_{1}$ and has inertia group $\bold
 Z_2$. 
 Thus the ramification $R$ of $p$ and $\omega_{\overline X}$ satisfy:
 $$\displaylines{4R =p^*( 2D_{1} + 3(D_{2}+D_{3}))\sim p^*(4L_2 +
   2D_{1} + D_2 + D_{3}) \cr
 2R \sim p^*(D_{1} + 3L_2) \sim p^*(2L_1 + 2L_2 - D_{2}) \sim
 p^*(2L_3 + D_{2}) \cr
 2K_{\overline X} \sim p^*(2K_Y + 2L_1 + 2L_2 - D_{2}) \sim p^*(2K_Y +
 2L_3 + D_{2})  \ \singleptone. \cr}$$

 If $D_2=0$, we have in fact
 that $\omega_{\overline X}=p^*(\omega_Y \otimes L_3)$ and $\overline
 X$ is Gorenstein. If $D_2=C_0$, then 
  $\omega_{\overline X}^{\otimes 2}=p^*(\omega_Y^{\otimes 2} \otimes
   L_3^{\otimes 2}
 \otimes \Cal O_Y(C_0))$ and $\overline X$ is $2$-Gorenstein.
 Let $F$ be the reduced cycle 
 consisting of the curves of $\overline X$
 lying over $C_0$. 
 Then $F$ is the exceptional locus of $\overline q$
 and according to formula \prsinglemmatwo \ns, 
 there exist suitable canonical divisors
 $K_X$ and $K_{\overline X}$ such that
 $K_{\overline X} \equiv {\overline q}^*K_X + aF$, with $a$ a nonnegative
 rational number.
 Since $\overline X$ is at worst locally
 $2$-Gorenstein, we have
 $2K_{\overline X} = {\overline q}^*2K_X + 2aF$, with $2a$ a 
nonnegative
 integer. Now we will determine $a$ and prove that $W=S(0,2)$ if $a >
 0$. 
We split the remaining of Case 1 into four subcases, according to
whether $C_0$ is in the branch locus of $p$ or not and 
according to what $D_2$ is:
 
\smallskip

{\it Case 1.1:} $C_0$ is not in the
 branch locus of $p$.  In this case $D_2=0$ and, as previously observed,
$L_1
 \otimes L_2 = L_3$ and $\overline X$ is Gorenstein. 
 Then $K_{\overline X} = {\overline q}^*K_X + aF$ and $a$ is a
 nonnegative integer. Moreover 
 $F=p^*C_0$ and is therefore a reduced Cartier divisor
 such that $F^2=-4e$. By \ratlemma we know that the arithmetic genus of
 $F$ is $0$, hence by adjunction
 $$-2=((K_{\overline X}+F) \cdot F)=({\overline q}^*K_X + (a+1)F)\cdot
 F)=-4e(a+1) \ .$$
 This gives $1=2e(a+1)$ but this is not possible since $a$ and $e$ are
 integers. Thus Case 1.1 does not occur.
 
\smallskip
 
{\it Case 1.2:} $C_0$ is in the branch locus of $p$ and
 $D_2=C_0$. 
Then the inertia group of $F$ is $\bold Z_4$, 
 $p^*C_0=4F$ and $F^2=-\frac e 4$. 
 By the previous observation
  $\omega_{\overline X}^{\otimes 2}$ is the pullback
 by $p$ of a certain line
 bundle on $Y$. On the other hand recall that $2K_{\overline X}$ is
 equal to ${\overline q}^*2K_X + 2aF$ and, since 
 $\omega_X=\varphi^*\Cal O_W(1)$ and by the
 commutativity of diagram \text{\diagramone \ns,} linearly equivalent to  
 $p^*(2C_0+2ef) + 2aF$. Thus $\Cal O_{\overline X}(2aF)=p^*N$, for certain
 line bundle $N$ on $Y$. Then $p^*(N^{\otimes 4})=p^*\Cal
 O_Y(2aC_0)$. This implies that $N^{\otimes 4}$ and $\Cal
 O_Y(2aC_0)$ are numerically equivalent, and since $Y$ is a rational
 ruled surface, linearly equivalent. Then $N=\Cal O_Y(\alpha C_0)$,
 with $4\alpha=2a$, and $\alpha$ integer so $a$ is in fact 
 a nonnegative even integer.
 On the other hand we consider the ramification $R$ of $p$. We have
 $$K_{\overline X} = p^*K_Y + R=p^*(-2C_0-(e+2)f)+R \ .$$
 Since $K_{\overline X} \equiv p^*(C_0+ef)+ aF$, we obtain $R \equiv
 p^*(3C_0+(2e+2)f)+aF$. Recall that $C_0$ is in the branch locus of
 $p$ and since neither $D_1+D_2$ nor $D_1+D_3$ nor $D_2+D_3$ has
 multiple components, $C_0$ belongs to the branch
 locus with multiplicity $1$. Then we can write
 $R=R_1+3F$, where $R_1$ is a cycle that does not contain $F$ in its
 support, and therefore, $R_1 \cdot F \geq 0$. Now we compute exactly
 $R_1 \cdot F$. The cycle $R_1$ is numerically equivalent to
 $p^*(2C_0+(2e+2)f)+(a+1)F$, then 
 
 $$R_1 \cdot F= (2C_0+(2e+2)f)C_0+(a+1)F^2=2-\frac{(a+1)e}{4} \ .$$
 Hence $(a+1)e \leq 8$. Then, since $a \geq 0$ and is even, 
 and $e \geq 2$, then either $a=0$ and $\overline q$ is crepant or
 $a=e=2$. The first case is not possible, since then $\omega_{\overline
   X}=p^*\Cal O_Y(1)$ and $L_1 \otimes L_2 =
 L_3$, hence $D_{2}=0$. In the second case 
we know that
$\overline X$ is
 $2$-Gorenstein. Moreover, since $e=2$, $F^2=-\frac e 4=-\frac 1 2$. 
 
\smallskip

 {\it Case 1.3:} $C_0$ is in the branch locus of $p$, $D_2=0$ and $F$ has
inertia group $\bold Z_4$ (last condition occurs if and only if $C_0 \subset
 D_{3}$). In this
 case $\overline X$ is locally Gorenstein.
We have $p^*C_0=4F$ and $F^2=-\frac e 4$. 
 By the previous observations 
 $\omega_{\overline X}$ is the pullback
 by $p$ of a certain line
 bundle on $Y$. On the other hand recall that $K_{\overline X}$ is
 equal to ${\overline q}^*K_X + aF$ (with $a$ nonnegative integer) 
and, since 
 $\omega_X=\varphi^*\Cal O_W(1)$ and by the
 commutativity of diagram \text{\diagramone \ns,} linearly equivalent
to  
 $p^*(C_0+ef) + aF$. Thus $\Cal O_{\overline X}(aF)=p^*N$, for certain
 line bundle $N$ on $Y$. Then $p^*(N^{\otimes 4})=p^*\Cal
 O_Y(aC_0)$. This implies that $N^{\otimes 4}$ and $\Cal
 O_Y(aC_0)$ are numerically equivalent, and since $Y$ is a rational
 ruled surface, linearly equivalent. Then $N=\Cal O_Y(\alpha C_0)$,
 with $4\alpha=a$, so $a$ is multiple of $4$.
 On the other hand we consider the ramification $R$ of $p$. We have
 $$K_{\overline X} = p^*K_Y + R=p^*(-2C_0-(e+2)f)+R \ .$$
 Since $K_{\overline X} \sim  p^*(C_0+ef)+ aF$, we obtain $R \sim
 p^*(3C_0+(2e+2)f)+aF$. Recall that $C_0$ is in the branch locus of
 $p$ and as argued before, $C_0$ belongs to the branch
 locus with multiplicity $1$. Then we can write
 $R=R_1+3F$, where $R_1$ is a cycle that does not contain $F$ in its
 support, and therefore, $R_1 \cdot F \geq 0$. Now we compute exactly
 $R_1 \cdot F$. The cycle $R_1$ is linearly equivalent to
 $p^*(2C_0+(2e+2)f)+(a+1)F$, then 
 
 $$R_1 \cdot F = (2C_0+(2e+2)f)C_0+(a+1)F^2=2-\frac{(a+1)e}{4} \ .$$
 Hence $(a+1)e \leq 8$. Then, since $a \geq 0$ and is multiple of $4$
 and $e \geq 2$, then we have $a=0$. 
 
\smallskip

{\it Case 1.4:} $C_0$ is in the branch locus of $p$, $D_2=0$ and $F$ has
inertia group $\bold Z_2$ (last condition holds if and only if $C_0
\subset D_{1}$).   
 In this
 case $\overline X$ is also Gorenstein.
We have $2F=p^*C_0$  and $F^2=- e$. 
 By the previous observations 
 $\omega_{\overline X}$ is the pullback
 by $p$ of a certain line
 bundle on $Y$. On the other hand recall that $K_{\overline X}$ is
 equal to ${\overline q}^*K_X + aF$ (with $a$ a nonnegative integer) and, 
since $\omega_X=\varphi^*\Cal O_W(1)$ and by the
 commutativity of diagram \text{\diagramone \ns,} linearly equivalent
to  
 $p^*(C_0+ef) + aF$. Thus $\Cal O_{\overline X}(aF)=p^*N$, for certain
 line bundle $N$ on $Y$. Then $p^*(N^{\otimes 2})=p^*\Cal
 O_Y(aC_0)$. This implies that $N^{\otimes 2}$ and $\Cal
 O_Y(aC_0)$ are numerically equivalent, and since $Y$ is a rational
 ruled surface, linearly equivalent. Then $N=\Cal O_Y(\alpha C_0)$,
 with $2\alpha=a$, so $a$ is even. 
 On the other hand we consider the ramification $R$ of $p$. We have
 $$K_{\overline X} = p^*K_Y + R=p^*(-2C_0-(e+2)f)+R \ .$$
 Since $K_{\overline X} \sim  p^*(C_0+ef)+ aF$, we obtain $R \sim
 p^*(3C_0+(2e+2)f)+aF$. Recall that $C_0$ is in the branch locus of
 $p$ and, as argued before, $C_0$ belongs to the branch
 locus with multiplicity $1$. Then we can write
 $R=R_1+F$, where $R_1$ is a cycle that does not contain $F$ in its
 support, and therefore, $R_1 \cdot F \geq 0$. Now we compute exactly
 $R_1 \cdot F$. The cycle $R_1$ is linearly equivalent to
 $p^*(2C_0+(2e+2)f)+(a+1)F$, then 
 
 $$R_1 \cdot F = 2(2C_0+(2e+2)f)C_0+(a+1)F^2=4-{(a+1)e} \ .$$
 Hence $(a+1)e \leq 4$. Since $a \geq 0$ and even
 and $e \geq 2$, then we have $a=0$. This ends our
argument when $G=\bold Z_4$.

 \medskip 
 
 {\it Case 2:} $G=\bold Z_2 \times \bold Z_2$. 
 
 Let $G_1$, $G_2$ and $G_3$ be the three order $2$ subgroups of $G$ and 
let $X_i$
 be the quotient of $X$ by $G_i$. As we argued in Case 2 of the proof
of \singlemma \ns, 
 associated to each subgroup $G_i$ we
 have a way of decomposing $\varphi$, namely, 
 $$X @> \varphi_2^i >> X_i @> \varphi_1^i >> W \ ,$$ 
 where $X_i$ is normal and $\varphi_1^i$ and $\varphi_2^i$ are Galois
 covers with Galois group $G/G_i$ and $G_i$ respectively.
 
 Let $\overline{X_i}$ be such that $\Cal O_{\overline X_i}$ is the integral
 closure of $\Cal O_Y$ in $\Cal K(X_i)$. Then, arguing as in Case 2 of
the proof of \singlemma we
 have the following commutative diagram
 
 $$ \matrix \overline X & @> \overline q >> & X \cr
 @VV p_2^i V &  \hskip -.8 truecm @VV\varphi_2^i V \cr
 \overline {X_i} & @> q_i >> & X_i \cr
         @VV p_1^i V &  \hskip -.8 truecm @VV\varphi_1^i V \cr
 Y & @> q >>& W \cr 
 \endmatrix $$  
 where $p_1^i$ and $p_2^i$ are Galois covers with groups $G/G_i$ and
 $G_i$ respectively and $p=p_1^i \circ p_2^i$. Moreover, by
 construction, $\overline{X_i}$ is normal, and, since $Y$ is
 smooth and $p_1^i$ is a double cover, it is also locally Gorenstein. 
 
\smallskip

 We will show now that $W=S(0,2)$.
 We examine the structure of the branch locus of $p$ and of $p_1^1$,
 $p_1^2$ and $p_1^3$. Since $Y$ is smooth, $p$ is a flat Galois cover with
 group $\bold Z_2 \times \bold Z_2$. Recall that by \Galquad there exist
 divisors $D_1$, $D_2$ and $D_3$ in $Y$ such that $p_1^i$ is a double
 cover of $Y$ branched along $D_j+D_k$, where
 $i \neq j, j \neq k$ and $k \neq i$. 
 Then since $\overline{X_i}$ is normal for all
 $i=1,2,3$, no two among $D_1$, $D_2$ and $D_3$ have a common
 component. In particular $C_0$ is not contained in two of them, let us
 say, $D_2$ and $D_3$, so $C_0$ is not in the branch locus of $p_1^1$. 
 Let $E={p_1^1}^*C_0$. Then $E$ is a reduced Cartier divisor and the
 exceptional locus of $q_1$ and $E^2=-2e$. Now,
 deg$\varphi_2^i(\text{disc}X+1) \geq \text{disc}X_i + 1$ (see \CKM \ns,
 6.7.i; note the statement in \CKM is incorrect: the morphism should
 be required to be finite). Now, since $X$ has canonical
 singularities, the discrepancy of each of the 
 $X_i$ is greater than or equal to $- \frac 1 2$, so in particular
 $X_i$ has rational singularities.  Then by \ratlemma \ns, 
 $E$ has arithmetic genus
 $0$. On the other hand since $\overline{X_1}$ is normal and
 locally Gorenstein, for suitable canonical divisors $K_{\overline{X_1}}$ and
 $K_{X_1}$, we obtain as in \prsinglemmaone 
 $$K_{\overline{X_1}}\equiv q_1^*K_{{X_1}} + bE \ ,$$
 where $b$ is a nonnegative rational number. Using adjunction we obtain
 $$-2=(K_{\overline{X_1}}+E)\cdot E = (q_1^*K_{{X_1}} +
 (b+1)E) \cdot E)= -2(b+1)e \ ,$$
 hence $b=-\frac{e-1}{e}$. Resolving the singularities of
 $\overline{X_1}$ and composing with $q_1$ we obtain a resolution of
 singularities for $X_1$. Since disc$X_1 \geq -\frac 1 2$, then $b \geq
 - \frac 1 2$ and $e=2$. 

\smallskip

Now we see that $\overline q$ is crepant. Let
$F$ be the reduced cycle consisting 
 of the curves of $\overline X$ lying over $C_0$. 
 We will show now that $C_0$ is in the
 branch locus of $p$. Assume it is not. Then 
 $F=p^*C_0$ and, by formula \prsinglemmaone \ns, 
 for suitable canonical divisors $K_X$ and $K_{\overline
   X}$ we have $$K_{\overline
   X} \equiv {\overline q}^*K_X+aF=p^*(C_0+2F) \ ,$$ with $a$
 a nonnegative rational 
 number and 
 $$K_{\overline
   X}=p^*K_Y+R \sim p^*(-2C_0-4f)+R \ ,$$
 where $R$ is the ramification divisor of $p$. 
 Then $R \equiv p^*(3C_0+6f) + aF$.
 Now if we denote $E_i={p_1^i}^*C_0$, by Zariski's main theorem, $E_i$ is
 connected, and, as pointed out before, $E_i$ is reduced and by \ratlemma \ns, 
 has arithmetic genus $0$. 
 Then $$(D_j+D_k) \cdot C_0=2 \ \prsinglemmaseven $$ 
 and
 since $2R=p^*(D_1+D_2+D_3)$, $R \cdot F =6$. But then
 $$6 = R \cdot F = p^*(3C_0+6f) \cdot p^*C_0 + a(p^*C_0)^2=-8a \ ,$$
 hence $a=-\frac 3 4$, which contradicts the fact that $a \geq 0$. 
 Hence we have seen that $C_0$ is in the branch locus of $p$. Then
 $C_0$ is contained in one of the $D_i$s, let us say, $D_1$, and, since
the $\overline X_i$s are normal, 
 $C_0 \not\subset D_2$, $C_0 \not\subset D_3$. To abridge our notation we
 set $p_1=p_1^1$ and $p_2=p_2^1$. Now by the same argument as the
 one used 
 just before \prsinglemmaseven $E=p_1^*C_0$ is a reduced, connected 
 Cartier divisor of arithmetic genus $0$ and $(D_2+D_3) \cdot C_0=2$. 
 Then $E$ is in the branch locus of $p_2$, otherwise $C_0$ would not 
 be in the branch locus of $p$, and we denote by $F$ the inverse 
 image of $E$ by $p_2$. Then $2F=p_2^*E=p^*C_0$ and $F$ is reduced 
 and connected, and has one or two components, depending on
 whether $E$ has one or two components. In any case, by
 \prsinglemmaone \ns, 
 for suitable
 canonical divisors $K_X$ and $K_{\overline X}$ we have
 $$\displaylines{ 
 K_{\overline X}=p^*K_Y+R \sim p^*(-2C_0-4f) + R_1 + F \cr
 K_{\overline X} \equiv {\overline q}^*K_X +
aF \sim
   p^*(C_0+2f) + aF  \ ,  \quad \rlap \prsinglemmaeight
}$$
 where $a$ is a nonnegative rational number, $R$ is the ramification
 divisor of $p$ and $R=R_1+F$ with $F$ not in the support of $R_1$. 
 Since $p^*C_0=2F$, this yields 
 $$R_1 \equiv p^*(2C_0+6f) + (a+1)F \ . \prsinglemmanine $$
 Now $(D_2 + D_3) \cdot C_0=2$
 implies $R_1 \cdot F \geq 2$. Then 
 $$\displaylines{2 \leq R_1 \cdot F = (p^*(2C_0+6f) +
   \frac{a+1}{2}p^*C_0) \cdot 
 \frac 1 2 p^*C_0  
 = \cr 2(2C_0+6f)\cdot C_0 + (a+1)C_0^2 = 4 -2(a+1) \ .}$$
 This implies $a \leq 0$, therefore $a=0$, i.e., $\overline q$ is crepant. 

\smallskip

Finally we prove that $\overline X$ has canonical singularities. If we
compose a
 resolution of singularities of $\overline X$ with $\overline q$ we
 obtain a resolution of singularities of $X$. Since $a=0$, 
 \prsinglemmaeight becomes
 $K_{\overline X} \equiv {\overline q}^*K_X $ and the discrepancies of
 the exceptional divisors of the resolution of $\overline X$ are the
 same whether considered with respect to $X$ or with respect to
 $\overline X$. Since $X$ has canonical singularities, disc$X \geq 0$
 and so disc$\overline X \geq 0$. Thus $\overline X$ has also canonical
 singularities.
 \qed  
 
\heading 3. Quadruple Galois canonical covers: crepant case \endheading

In this and the next section we achieve the classification of
quadruple Galois canonical covers $\varphi$ of singular rational
normal scrolls $W$. 
In the previous section we constructed a desingularization diagram for
$\varphi$ (see \diagramone \ns) and,
in \singlemma and \thmsinglept \ns, we
studied in great detail one of the sides of this diagram, namely, the
morphism $\overline q$. In \singlemma and \thmsinglept we came to the
following conclusion: either $\overline q$ is crepant, that is,
$\overline q^*\omega_X=\omega_{\overline X}$, or it is not, but in the
latter case, we do know many things about the discrepancies of
$\overline q$, about what the Galois group of $\varphi$ is and of what
kind the possible singularities of $\overline X$ are. Thus we will
split the study of $\varphi$ in two cases: the case in which
$\overline q$ is crepant and the case in which it is not. We deal with
the former case in this section and we will deal with the latter case
in Section 
4.  

\medskip

Let $\varphi, q, p$ and $\overline q$ and let 
$\overline \varphi = q \circ p = \varphi \circ \overline q$. 
If $\overline q$ is crepant, we can ``morally'' think of $\overline
\varphi$ as a Galois canonical cover. Admittedly, $\overline \varphi$
is not finite, so it is not a Galois cover according to our
definition, but since $\omega_{\overline X}=\overline q^*\omega_X$, it
turns out that $\omega_{\overline X}$ and, in fact, $\overline
\varphi$ is the canonical morphism of $\overline X$. However
$\omega_{\overline X}$ is not ample, so eventually $\overline \varphi$
maps $\overline X$ not onto an isomorphic image of the Hirzebruch
surface $Y$, but onto a ``singular realization'' of $Y$, namely, the 
singular rational normal scroll $W$. All this suggests that one can
deal with canonical covers $\varphi$ when $\overline q$ in way
parallel to the study of Galois canonical covers of smooth rational
normal scrolls carried out in \GPsmooth \ns. To do so we start giving
a definition:   

\proclaim{\convention}  Let $\overline X$ be a normal surface of general
type with canonical
singularities whose canonical line bundle is
base-point-free and let  $\overline X @> \overline \varphi >> W$ be the
canonical morphism of $\overline X$. We say that $\overline  \varphi$ 
satisfies \convnumber if it factorizes as follows:
 $$\overline X @> p >> Y @> q >> W \ ,$$
 where $p$ is finite  and 
 $q$ is the 
 minimal desingularization of $W$.
 \endproclaim 
 
In \GPsmooth we proved some general results concerning finite canonical
covers and Galois canonical covers  
of smooth surfaces of minimal degree. This results hold in
slightly greater generality, as they hold both for the above mentioned
canonical
covers of smooth surfaces and for finite covers $p$ as in \convention
\ns. Thus we proceed now to state the versions of the required results 
of \GPsmooth for morphisms $p$ and $\overline \varphi$ like those in
\convention \ns:

\proclaim{\split} Let $\overline X$ be a normal surface of general
type with canonical
singularities and base-point-free canonical bundle. 
Assume that the canonical morphism $\overline X @> \overline \varphi
>> W$ 
satisfies \convnumber \ns.  
Let
$H=q^*\Cal O_W(1)$.

\item{1)}  If $p$ has degree $4$, 
then $p_*\Cal O_{\overline X}$ is a vector bundle on $Y$ and 
 $$p_*\Cal O_{\overline X}=\Cal O_Y \oplus E \oplus (\omega_Y
 \otimes H^*)$$
  with $E$ a vector bundle over $Y$ of rank $2$. 

\item{2)} 
 If, in addition to the hypothesis in 1), $p_*\Cal O_{\overline X}$
 splits as a sum of line
 bundles, then 
 $$p_*\Cal O_{\overline X}=\Cal O_Y \oplus L_1^* \oplus L_2^* \oplus (\omega_Y
 \otimes H^*) $$ 
 with $L_1^* \otimes L_2^*= \omega_Y
 \otimes H^*$. 

\item{3)} If, in addition  to the hypothesis in 1), $\overline X$ is regular,  
then $$\displaylines{p_*\Cal O_{\overline X}= \Cal O_Y \oplus \Cal
 O_Y(-C_0-(e+1)f) \oplus \cr \Cal O_Y(-2C_0-(e+1)f) \oplus  
 \Cal O_Y(-3C_0-(2e+2)f) \ .}$$   
 \endproclaim
 
 \noindent {\it Sketch of proof.} This proposition is analogous to
 \GPsmooth \ns, Propositions 1.3 and 1.6, 3) and one can go through
 the proofs there and adapt them to the present situation. 
We make explicit below
 the parallelism between the two settings. 
In our case, $Y$, which  is a smooth Hirzebruch surface $\bold F_e$
 with $e \geq 
 2$, plays the role of $W$ in \GPsmooth \ns. The morphism $p$ (which
 is flat and finite, so  
 $p_*\Cal O_{\overline X}$ is a vector bundle over $\Cal O_Y$ of rank
 $4$) plays the role of $\varphi$ in \GPsmooth. Finally the role
 played in \GPsmooth by 
 the line bundle $\Cal O_W(1)=\Cal O_W(C_0+mf)$ is here played
 by the line 
 bundle $H=q^*\Cal O_W(1)=\Cal O_Y(C_0+ef)$. Since $\overline \varphi$
 is the canonical morphism of $\overline X$, the canonical bundle of
 $\overline X$ is $\omega_{\overline X}=p^*H$, and we can use relative
 duality for $p$ as we did for $\varphi$ in \GPsmooth \ns. 
The fact that $\overline \varphi$ is induced by the complete canonical
 series implies that  $H^0(p^*\Cal O_Y(C_0+ef))=H^0(\Cal
 O_Y(C_0+ef))$. The regularity of $\overline X$ assumed in 3) has the
 same implications for the summands of $p_*\Cal O_{\overline X}$ as
 the regularity of $X$
 has for the summands of $\varphi_*\Cal O_X$ in \GPsmooth \ns. 
\qed

\proclaim{\biprocyc }  Let $\overline X$ be a normal surface of general
type with canonical
singularities and base-point-free canonical bundle. 
Assume that the canonical morphism $\overline X @> \overline \varphi
>> W$ 
satisfies \convnumber  and that $p$ is Galois with Galois group
$G$. Let $L_1$ and $L_2$ be as in \split \ns, 2). 

\item{1)} If $G=
\bold Z_2
 \times \bold Z_2$, 
 then $\overline X$ is the fiber product over $Y$ 
 of two double covers $\overline X_1 @> p_1 >> Y$ and $\overline X_2
 @> p_2 >> Y$ 
 and $p$ is the natural map from the
 fiber product to $Y$. The trace-zero modules of $p_1$ and $p_2$ are
 $L_1$ and $L_2$.

\item{2)} If $G = \bold Z_4$, then there are two divisors $D_1$ and
  $D_2$ on $Y$ such that $p$ is the composition of a flat
  double cover $\overline X_1 @> p_1 >> 
 Y$ branched along $D_{2}$ followed by a flat double cover
 $\overline X 
 @> p_2 >> \overline X_1$, 
 branched along $p_1^*D_{1}$ and the ramification
 locus of $p_1$. Moreover, the trace zero module of $p_2$ is
 $p_1^*L_1$ and the trace zero module of $p_1$ is $L_2$. 

\endproclaim
 
 {\it Sketch of proof.} This result is analogous to  \GPsmooth \ns,
 Proposition 2.6, 2)  and 2.7, 4). In our setting, $\overline X$ and
 $Y$ are normal varieties and $p$ is a flat,
 Galois cover, so \Galquad applies to $p$ as it does to $\varphi$ in
 \GPsmooth \ns. 
Since $\overline \varphi$ satisfies \convnumb \ns, then \split implies
 that there is a splitting 
$$p_*\Cal O_{\overline X}=\Cal O_Y \oplus L_1^* \oplus L_2^* \oplus
L_3^* \ ,$$
with $L_1 \otimes L_2 =
 L_3= 
 \omega_Y^*
 \otimes H$. Then \GPsmooth \ns, Propositions 2.6, 1), 2.7, 1) and 2)
 apply to $\overline X @> p >> Y$. \qed

\medskip

Now we are ready to classify of quadruple Galois
canonical covers $\varphi$ of singular rational normal scrolls $W$
when the morphism $\overline q$ defined in \diagram is crepant.  
To each cover $X @> \varphi >> W$ there corresponds a
unique cover $\overline X @> p >> Y$ and we will classify these
latter covers. 
We will study separately the cyclic and the bidouble case:

 \proclaim{\crepantbid}
 Let  $W$ be a
 singular rational normal scroll and let $X 
 @> \varphi >> W$ be a Galois canonical cover with Galois
 group $\bold Z_2 \times \bold Z_2$.  Let $\overline X$, $Y$, $q$,
 $\overline q$ and
 $p$ be as in 
 \diagramone \ns. 
If $\overline q$ is crepant, then
\item{1)} $W=S(0,2)$ (and hence, 
 $Y=\bold F_2$);
\item{2)} $\overline X$ has at
 worst canonical 
 singularities; 
\item{3)} $X$ is regular;
 \item{4)} $\overline X @> \overline q >> X$ is the morphism from
   $\overline X$
to its canonical model; 
 \item{5)} $\overline X$ is the fiber product over $Y$ of two double 
 covers $p_1$ and $p_2$ branched along divisors
 $D_2$ and $D_1$ which are linearly equivalent to $2C_0+6f$ and
 $4C_0+6f$ respectively.

 \smallskip

\noindent  Conversely, let $\overline X$ be a normal surface with at
worst
 canonical singularities
  and let $Y=\bold
 F_2$. If
 $\overline X @> p >> Y$ is a fiber product of  two double 
 covers $p_1$ and $p_2$ as described in 5) above,
then 
  there exists a commutative diagram like \diagramone 
where $\overline q$
 is crepant and $\varphi$ is the canonical morphism of $X$ and 
 is Galois with Galois group $\bold Z_2
 \times \bold Z_2$. 
 \endproclaim
 
 \noindent {\it Proof.} It follows from \singlemma and \thmsinglept
that, 
 if $\overline q$ is crepant,
 then $\overline X$ has at worst canonical singularities, so we have 2).

\smallskip

Since $\overline q$ is
crepant, $H^0(\omega_{\overline X})=H^0(\omega_X)$ and
$\omega_{\overline X}$ is base-point-free, so
$\overline
\varphi=\varphi \circ \overline q$ is the canonical
morphism of $\overline X$. Then $\overline X
 @> \overline 
 \varphi >> W$, $p$ and $q$  satisfy the hypothesis of  \convention and  
  $Y=\bold F_e$.
The morphism $\varphi$ is Galois
 and  so is $p$, by  \samegrplem  and both have the same Galois
group, also by \samegrplem \ns. 
 Applying \bipro
 to $\overline X @> \overline
 \varphi >> W$, we obtain that $\overline X$ is the
 fiber product over $Y$  of two double covers $\overline X_1 @> p_1
 >> Y$ and
 $\overline X_2 @> p_2 >> Y$. Those covers are branched along divisors 
$D_1$ and $D_2$. We set $L_1$ and $L_2$ 
 line bundles on $Y$ such that $L_i=\Cal O_Y(a_iC_0+b_if)$,
 $L_1^{\otimes 2}=\Cal 
O_Y(D_2)$ and $L_2^{\otimes 2}=\Cal O_Y(D_1)$ (we are using this
 rather strange notation so as to be consistent with the notation of
 \Galquad \ns). Then 
 $$p_*\Cal O_{\overline X}= \Cal O_Y \oplus L_1^* \oplus L_2^* \oplus
 (L_1^* \otimes L_2^*) \ .$$

\smallskip

We prove now 1) and the remaining of 5), that is, the description of
$D_1$ and $D_2$ (we have already seen above that $\overline X$ is a
fiber 
product of two covers).    
By \biprocyc \ns, 1) and \split \ns, 2), we know that
$L_1^* \otimes L_2^*=\omega_Y \otimes H^*$, hence
$$\displaylines{a_1 + a_2 = 3 \cr
 b_1 + b_2 = 2e + 2 \ .} $$
 Since $D_1$ and $D_2$ are effective and linearly equivalent 
to $2(a_2C_0+b_2f)$ and $2(a_1C_0+b_1f)$ respectively,
 we have $a_i, b_i \geq 0$. We set $a_1=0,1$ (in which case,
 $a_2=3,2$). Since $\overline \varphi$ is induced by the complete
 canonical series of $\overline X$, $\overline q$ is crepant,
 $\varphi$ is a canonical cover and $\overline \varphi=q \circ p =
 \varphi \circ \overline q$, then $H^0(p^*\Cal O_Y(C_0+ef))=H^0(\Cal
 O_Y(C_0+ef))$. Then $H^0(\Cal O_Y((1-a_1)C_0+(e-b_1)f)=0$, so $b_1
 \geq 
 e +1 $ and, since $ b_1 + b_2 = 2e + 2$, $b_2 \leq e+1$. Now, let us
 assume $a_1=0$. Then $D_1 \sim 6C_0+2b_2f$, and, since $e \geq 2$,
 then $3C_0$ is in the fixed part of $|D_2|$. This would imply that
 $\overline X$ is nonnormal, which is not possible. Thus $a_1$ can
 only be $1$. Then $D_1 \sim 4C_0 + 2b_2f$.  If $e > 2$ or $b_2 <
 e+1$, then  $2C_0$ 
 is in the fixed part of $|D_2|$, and as before this is not possible. 
Thus we conclude that $a_1=1, a_2=2, b_1=b_2=e+1$ and $e=2$, and since
$L_i^{\otimes 2}=\Cal O_Y(D_j)$, we get $D_2 \sim
2C_0+6f$ and $D_1 \sim 4C_0+6f$.

\smallskip

Now we prove 3), that is, we show that $X$ is regular. The
irregularity of $\overline X$ is the sum of
$h^1(\Cal O_Y)$, $h^1(L_1^*)$, $h^1(L_2^*)$ and $h^1(L_1^* \otimes
L_2^*)$, and those numbers are $0$ for the above values of $a_1, a_2,
b_1$ and $b_2$.  Therefore $\overline X$ is regular and, since
$\overline X$ and $X$ are birational and have rational singularities,
so is $X$. 

\smallskip

  Finally we show 4). 
Recall that
 $\overline \varphi$ is the canonical morphism of $\overline X$. Since
 $\varphi$ is finite, the curves (which are $-2$-curves) contracted by
 $\overline \varphi$ are the same as the curves contracted by
 $\overline q$. Thus $\overline X @> \overline q >> X$ is the morphism
from $\overline X$ to its 
canonical model.

 \smallskip
 
 Conversely, assume that $\overline X$ is a normal surface with
 canonical singularities and that
 $\overline X @> p >> Y$ is the
 fiber product over $Y=\bold F_2$ of two double covers of $Y$, 
one branched along
 a divisor linearly equivalent to $4C_0+6f$ and the other
 branched along a divisor linearly equivalent to $2C_0+6f$.
\biproconv tells that 
$p$ is a
 Galois cover with group $\bold Z_2 \times \bold Z_2$. \biproconv
also tells that  the canonical bundle of $\overline X$ is $p^*\Cal
O_Y(C_0+2f)$, so it is
 base-point-free. Moreover one sees easily using projection formula that
 $H^0(\omega_{\overline X})=H^0(\Cal O_Y(C_0+2f))$, so 
 the canonical morphism $\overline \varphi$ of $\overline X$
  factors as $\overline \varphi = q \circ p$,
 where $Y @> q >> W$ is the contraction of $C_0$. 
 Now let $X$ be the
 canonical model of $\overline X$. 
 Then $\overline \varphi$ also
 factors as $\overline \varphi = \varphi \circ \overline q$, where $X @
 > \varphi >> W$ is the canonical morphism of $\overline X$. 
 Finally, we note that, since the canonical bundle of $\overline X$ is
 base-point-free, $\overline q$ is crepant, and, since $p$ is Galois
 with group $\bold Z_2 \times \bold Z_2$, so is $\varphi$, by
 \samegrplem \ns.  
 \qed

 \proclaim{\crepantcyc}
 Let $W$ be a
 singular rational normal scroll and let $X 
 @> \varphi >> W$ be a Galois canonical cover with Galois
 group $\bold Z_4$. 
 Let $\overline X$, $Y$, $q$, $\overline q$ and $p$ be as in 
 \diagram \ns. 
If $\overline q$ is crepant, then

\item{1)} $W=S(0,2)$ (and hence, 
 $Y=\bold F_2$); 

\item{2)} $\overline X$ has at
 worst canonical 
 singularities; 

\item{3)} $X$ is regular;

\item{4)}  $\overline
 X @> \overline q >> X$ is the morphism from $\overline X$ to its
canonical model;

 \item{5)} $p$ is the composition of two
 double covers  
 $\overline X_1 @> p_1 >>Y$ branched along a divisor $D_2$ linearly
 equivalent to $4C_0+6f$ and $\overline X @> p_2 >> \overline X_1$,
 branched along the ramification of $p_1$ and $p_1^*D_1$, with $D_1$
 linearly equivalent to $3f$ and having trace-zero module $p_1^*\Cal
 O_Y(-C_0-3f)$. 
 
 \smallskip

 \medskip

 Conversely, let $\overline X$ be a normal surface with at worst
 canonical singularities and let $Y=\bold F_2$. If $\overline X @> p >>
 Y$ is the composition of two
 double covers  $p_1$ and $p_2$ as described in 5) above,
 then  there exists a commutative diagram like
               \diagramone  
 where $\overline q$
 is crepant and $\varphi$ is the canonical morphism of $X$ and 
 is Galois with Galois group $\bold Z_4$.

\endproclaim
 
 \noindent {\it Proof.} 
 It follows from \singlemma and \thmsinglept that, 
 if $\overline q$ is crepant,
 then $\overline X$ has at worst canonical singularities, so we have
2). 
 Let $\overline \varphi=\varphi \circ \overline q$. Then, as argued in
the proof of \crepantbid \ns, 
$\overline X
 @> \overline 
 \varphi >> W$ is a canonical cover satisfying \convnumber \ns.
   Let $Y, p$ and $q$
 satisfy the hypothesis of \text{\convention \ns.}
 The morphism $\varphi$ is Galois
 and, by \samegrplem \ns,  so is $p$ and both have the same Galois
 group, which is  $\bold 
 Z_4$. 
Applying \split \ns, 2) and \biprocyc \ns, 2) 
 to $\overline X @> \overline
 \varphi >> W$, we obtain that $\overline X$ is the
 composition of a flat
  double cover $\overline X_1 @> p_1 >> 
 Y$ branched along $D_{2}$ followed by a flat double cover
 $\overline X 
 @> p_2 >> \overline X_1$, 
 branched along $p_1^*D_{1}$ and the ramification
 locus of $p_1$. Moreover, the trace zero module of $p_2$ is
 $p_1^*L_1$ and the trace zero module of $p_1$ is $L_2$. 
Then, by \Galcycliconv \ns, 
$$p_*\Cal O_{\overline X}= \Cal O_Y \oplus L_1^* \oplus L_2^* \oplus
 (L_1^* \otimes L_2^*) \ .$$

\smallskip

We prove now 3), that is, that $X$ is regular.
  Let $L_i=\Cal O_Y(-a_iC_0-b_if)$. 
By \biprocyc \ns, 2) and \split \ns, 2), we know that
$L_1^* \otimes L_2^*=\omega_Y \otimes H^*$, hence
$$\displaylines{a_1 + a_2 = 3 \cr
 b_1 + b_2 = 2e + 2 \ .} $$
Let us assume that $\overline X$ is irregular.
Since $L_1^{\otimes 2} \otimes L_2^* = \Cal
 O_Y(D_1)$ and $L_2^{\otimes 2}=\Cal O_Y(D_2)$ are  effective then 
 $a_2 \leq 2a_1$, $b_2 \leq 2b_1$ and $a_2, b_2 \geq 0$. Then $a_1, b_1
 \geq 0$ also. Moreover, $a_1, b_1 \geq 1$, otherwise we will
 contradict $a_1 + a_2 = 3$ or $b_1 + b_2 = 2e + 2$. Then let us
 examine all the possibilities for $a_1$. First, if $a_1=1$, then
 $L_1^*=\Cal O_Y(-C_0-b_1f)$ and $L_3^*=\Cal O_Y(-3C_0-(2e+2)f)$ are
 special, so $L_2^*=\Cal O_Y(-2C_0-b_2f)$ is not. This means $b_2
 \leq e$. In this case, since $e \geq 2$, $2C_0$ is in the fixed part
 of $|D_2| = |4C_0 + 2b_2f|$, and $\overline X$ will not be normal,
 which is not possible, so $a_1=1$ is ruled out. Second, if $a_1=2$,  
 then $D_1 \sim 3C_0+(2b_1-b_2)f$, $D_2 \sim 2C_0+2b_2f$, $L_1 = \Cal
 O_W(2C_0+b_1f)$ and $L_2=\Cal
 O_W(C_0+b_2f)$. Since $H^1(L_2^*)=H^1(L_3^*)=0$ and we are assuming
 $\overline X$ to be irregular, then
 $H^1(L_1^*) \neq 0$. This implies $b_1 \leq e$. Then $b_1 + b_2 = 
 2e + 2$ implies $b_2 \geq e+2$. On the other hand, since $\overline X$
 is 
 normal, $C_0$ has at most multiplicity $1$ in the fixed part of
 $|3C_0+(2b_1-b_2)f|$, and this implies $2b_1-b_2-2e \geq 0$. Then we
 have $e \leq -2$, which is a contradiction, so $a_1=2$ is also ruled out. Then the only
 possibilities left is $a_1=3$.  If $a_1=3$, then  $D_1 \sim 6C_0 +
 (2b_1-b_2)f$, and since $X$ is normal,  
$C_0$ has at most multiplicity $1$ in the fixed part of $|D_1|$, hence 
 $2b_1 - b_2 -5e \geq 0$. Now since
 $H^0(p^*\Cal O_Y(C_0+ef))=H^0(\Cal O_Y(C_0+ef))$, 
we have that $b_2 >
 e$, hence $b_1 < e+2$. Then we get $-4e+4 > 0$. This contradicts $e
 \geq 2$ and so $a_1=3$ is also ruled out. Thus $\overline X$ is
 regular, and so is $X$. 
   
\smallskip

We prove now 1) and the remaining of 5), that is, the description of
$D_1$ and $D_2$. 
Since $\overline X$ is regular, then \split \ns, 3) tells that
$b_1=b_2=e+1$ and, either $a_1=1,
a_2=2$ or $a_1=2,
a_2=1$. If $a_1=2$ and $a_2=1$, 
then $L_2^{*}=\Cal
 O_Y(-C_0-(e+1)f), L_1^{*}=\Cal O_Y(-2C_0-(e+1)f)$ and 
 $D_{1}\sim 3C_0+(e+1)f$. Then, since $e \geq 2$, $2C_0$ is in the
 fixed part of $|D_1|$ and this contradicts the normality of
 $\overline X$. Therefore the only possibility left is $a_1=1,
a_2=2, b_1=b_2=e+1$. In this case $L_1^{*}=\Cal
 O_W(-C_0-(e+1)f), L_2^{*}=\Cal O_W(-2C_0-(e+1)f)$, so $D_{2}$ is
 linearly equivalent to $4C_0+(2e+2)f$. Arguing as before we see that
 the  normality of $\overline X$ implies $e \leq 2$, so in fact,
 $e=2$. Since in this case   $D_{1}$ is linearly equivalent to
 $(e+1)f$, this concludes the proof of 1) and 5).

\smallskip

 Finally, by the same argument given in the
 proof of \crepantbid \ns, $\overline X @> \overline q >> X$ is the
 morphism from $\overline X$ onto its canonical model $X$, so we
have 4).  
 
 \smallskip
 
 Conversely, let $Y=\bold F_2$, let $\overline X$ be a normal surface
 with at worst 
 canonical singularities and let 
 $\overline X @> p >> Y$ be the composition of two
 double covers $\overline X_1 @> p_1 >> Y$, branched along $D_2$ and
 $\overline X @> p_1 >> \overline X_1$, branched along $p_1^*D_1$ and
 the ramification divisor of $p_1$, where $D_1 \sim 3f$ and $D_2
 \sim 4C_0+6f$ and $p_1^*\Cal O_Y(-C_0-3f)$ is the trace
 zero module of $p_2$. 
By \Galcycliconv \ns,  $p$ is a Galois cover with group
 $\bold Z_4$, and by \samegrplem
\ns, so is $\varphi$. On the other hand, the trace zero module of
$p_2$ is $\Cal O_Y(-2C_0-3f)$, since $D_2 \sim 4C_0+6f$. Then
\Galcycliconv also implies that $\omega_{\overline X}= p^*\Cal
O_Y(C_0+2f)$. 
 Then $\omega_{\overline X}$ is base-point-free and one sees easily
 using projection formula that 
 $H^0(\omega_{\overline X})=H^0(\Cal O_Y(C_0+2f))$, so 
 the canonical morphism $\overline \varphi$
  factors as $\overline \varphi = q \circ p$,
 where $Y @> q >> W$ is the contraction of $C_0$. Now let $X$ be the
 canonical model of $\overline X$. 
 Then $\overline \varphi$ also
 factors as $\overline \varphi = \varphi \circ \overline q$, where $X @
 > \varphi >> W$ is the canonical morphism of $X$. 
 Since the canonical bundle of $\overline X$ is
 base-point-free $\overline q$ is crepant. 
 \qed
 
\heading 4. Quadruple Galois covers: non crepant case. \endheading

In this section we complete the classification of quadruple Galois
canonical covers $\varphi$ of singular rational normal scrolls $W$. Now we are
concerned with the case in which the morphism $\overline q$ in
\diagramone is non-crepant. As we did in Section 3, instead of
classifying directly canonical Galois covers $X @ > \varphi >> W$ we will
classify their ``desingularization'' $\overline X @> p >> Y$ (see
\diagram \ns) and we will
specify the way to go from $\overline X$ to $X$ and vice versa, by
explicitly characterizing the birational morphism $\overline q$.

\proclaim{\noncrepantbid} Let $W$ be a
 singular rational normal scroll and let $X 
 @> \varphi >> W$ be a Galois canonical cover with Galois
 group $\bold Z_2 \times \bold Z_2$. 
 Let $\overline X$, $Y$, $q$, $\overline q$ and $p$  be as in 
 \diagramone \ns.
If $\overline q$ is not crepant, 
 then

\item{1)} $W=S(0,2)$; 

\item{2)}  $\overline X$ has at worst canonical singularities;

\item{3)} $X$ is regular; 

\item{4)} 
$\overline X @>
\overline q >> X$ is the morphism from $\overline X$ to its canonical
model. 

\item{5)} $\overline
X$ is the
 normalization of the fiber 
 product over $Y$ of two double covers of $Y$ branched each along
 divisors $D_1'$ and $D_2'$, where $D_1'=D_1 + C_0$, $D_2'=D_2 +
 C_0$, $D_1 \sim D_2 \sim 3C_0+6f$  and all components of $D_1 +
 D_2 + C_0$ have multiplicity $1$.

 \medskip
 
 Conversely, if $\overline X$ has at worst canonical singularities and 
 is the
 normalization of a fiber 
 product over $Y$ as described in 5) above, 
 then there exists a
 commutative diagram like \diagramone such that $\overline q$ is
 noncrepant 
and
$\varphi$ 
 is the canonical morphism of $X$ and is Galois with Galois group $\bold
 Z_2 \times \bold Z_2$. 
 \endproclaim

 \noindent{\it Proof.} Since $\overline q$ is non crepant and $G = \bold Z_2
 \times \bold Z_2$, 1) and 2) follow from \singlemma and
 \thmsinglept \ns. We are in fact in the situation of \singlemma \ns,
 2), so 
$\overline q$ is the blowing down of two $-1$ curves, and, since 
$\omega_X$ is ample,
$X$ is not only minimal (in the sense that $K_X$ is nef) but is also its
canonical model. This shows 4).

Let us call
$F_1$ and $F_2$ the two curves (two smooth lines) lying over $C_0$. 
 Recall that, as seen in \singlemma just before \prsinglemmathree
\ns,  
 $p^*C_0=2F_1+2F_2$. 
 The ramification locus $R$ of $p$ is of the form $R_1 + F_1 + F_2$, 
 where $R_1$ contains neither $F_1$ nor $F_2$ and,  by \prsinglemmasixb
 \ns,  
 $R_1 \cdot F_i = 0$.
 By \prsinglemmatwo and having in account \singlemma \ns, 2), 
for suitable canonical
 divisors $K_{\overline X}$ and $K_X$ we have 
 $$K_{\overline X}={\overline q}^*K_X + F_1 +
 F_2  \quad \eqno \ncrone \ .  $$

 Let $G_1$, $G_2$ and $G_3$ be the three index $2$ subgroups of $G$ 
 and let $X_i$
 be the quotient of $X$ by $G_i$. As we argued in Case 2 of 
 \singlemma \ns, 
  associated to each subgroup $G_i$ we
 have a way of decomposing $\varphi$, namely, 
 $$X @> \varphi_2^i >> X_i @> \varphi_1^i >> W \ $$ 
 where $X_i$ is normal and $\varphi_1^i$ and $\varphi_2^i$ are Galois 
 covers of degree $2$. Let $\overline{X_i}$ be such that 
 $\Cal O_{\overline X_i}$ is the integral closure of $\Cal O_Y$ in 
 $\Cal K(X_i)$. Then, again arguing as in \singlemma  we
 have the following commutative diagram

 $$ \matrix \overline X & @> \overline q >> & X \cr
 @VV p_2^i V &  \hskip -.8 truecm @VV\varphi_2^i V \cr
 \overline {X_i} & @> q_i >> & X_i \cr
 @VV p_1^i V &  \hskip -.8 truecm @VV\varphi_1^i V \cr
 Y & @> q >>& W \cr 
 \endmatrix \  \ncrdiag $$  
 where $p_1^i$ and $p_2^i$ are Galois covers of degree $2$ and $p=p_1^i
 \circ p_2^i$. Moreover, by construction, $\overline{X_i}$ is normal
 and, since $Y$ is smooth and $p_1^i$ is a double cover, it is also 
 Gorenstein. Let 
 $$p_*\Cal O_{\overline X}=\Cal O_Y \oplus  L_1^* \oplus L_2^* \oplus
 L_3^* \ ,$$
 with $L_i=\Cal O_Y(a_iC_0+b_if)$. 
  Then by \Galquad \ns, 2) there exist effective divisors $D_1$, $D_2$ and
 $D_3$ on $Y$ such that $D_1+D_2+D_3$ is the branch locus of $p$,
 $\overline {X_i}=\text{Spec}(\Cal O_Y \oplus L_i^*)$, $L_i^{\otimes
 2}=\Cal O_Y(D_j+D_k)$  and
 $p_i^1$ is branched along
 $D_j+D_k$. By \ncrone and the commutativity of \ncrdiag we
 have
 $$\omega_{\overline X}=p^*\Cal O_Y(C_0+2f) \otimes \Cal
 O_{\overline X}(F_1+F_2) \ .$$
 
 We find out now the possible values of the $a_i$s and the $b_i$s. Let $D
\simeq
\bold P^1$ be a smooth general member in the linear
 system $|C_0+2f|$ on  
 $Y=F_2$ and let $C$ be
 its inverse image under $p$. Then, since $C_0 \cdot (C_0 + 2f) = 0$,
 by adjunction,  $\omega_C=p^*\Cal
 O_Y(2C_0+4f) \otimes \Cal O_C =p^*O_{\bold P^1}(4)$. Applying
 relative duality to $p|_C$ as we did 
 in Case 1 of the proof of \thmsinglept \ns, 
we conclude (maybe
renumbering $L_1$, $L_2$ and $L_3$) that 
 $(L_1 ^* \otimes L_2^*) \otimes \Cal O_D=L_3^* \otimes \Cal O_D=
 \Cal O_{\bold P^1}(-6)$. Then $b_3=6$ and $b_1+b_2=6$. 
 Then $b_1+b_2=b_3$ and, 
since $L_1 \otimes L_2 = L_3 \otimes \Cal O_Y(D_3)$, then $D_3 \cdot
D=0$, so $D_3$ is a multiple of $C_0$. Then  $D_3=0$ or
$D_3=C_0$, for $\overline X_1$ and $\overline X_2$ are normal. On
the other hand, if $D_3=0$, then $L_1 \otimes L_2 = L_3$, so by
\biproconv \ns, 
$\omega_{\overline X}=p^*\Cal O_Y(C_0+2f)$, a contradiction. Hence
we have $D_3=C_0$ and 
$$a_1 + a_2 = a_3 + 1 \quad \eqno \ncrbidthree \ .$$

Now,  since $\overline{X_i}$ 
 is normal, no two among $D_1$, $D_2$ and $D_3$ have common 
 components. Thus $2R=p^*(D_1+D_2+D_3)$.
 By \Galquad \ns, $(L_1 \otimes L_2 \otimes L_3)^{\otimes 2}=\Cal
 O_Y(2(D_1+D_2+D_3))$, and, since $Y$ is a rational ruled surface, 
 $L_1 \otimes L_2 \otimes L_3=\Cal
 O_Y(D_1+D_2+D_3)$. Thus from \prsinglemmatwo \ns, 
 \ncrone and the commutativity of \ncrdiag we
 have $$\displaylines{2K_{\overline X} \sim p^*(-4C_0-8f) + 2R \sim
 p^*((a_1+a_2+a_3-4)C_0+4f) \cr
 2K_{\overline X} \sim p^*(2C_0+4f) +2F_1 + 2F_2 = p^*(3C_0+4f) \ .}$$
 Then $(a_1+a_2+a_3-4)C_0+4f \equiv 3C_0+4f$, and since $Y$ is a
 rational normal scroll, $(a_1+a_2+a_3-4)C_0+4f = 3C_0+4f$, so
 $a_1+a_2+a_3=7$.  This together with \ncrbidthree yields
$a_1+a_2=4$ and $a_3=3$. We use this information, together with
$b_1+b_2=b_3=6$ previously obtained, to determine the $a_i$s and
the
$b_j$s. Since
$L_1^{\otimes 2} \otimes \Cal O_Y(-C_0)= \Cal O_Y(D_2)$ and
$L_2^{\otimes 2} \otimes \Cal O_Y(-C_0)= \Cal O_Y(D_1)$, and $D_1$
and $D_2$ are effective, we obtain $2a_1-1 \geq 0$ and $2a_2 -1 \geq
0$, yielding $a_1, a_2 \geq 1$ and also, $b_1, b_2 \geq 0$. Let us
assume $a_1 \leq a_2$. Then the only possibilities for $a_1, a_2$ are
$a_1=1$, $a_2=3$ or $a_1=a_2=2$. On the other hand $D_1$ and $D_2$ not
 having multiple components implies $b_1 \geq 2a_1 -1$ and $b_2 \geq 2a_2 -1$.
Thus
$a_1=1$, $a_2=3$ and $b_1+b_2=6$ implies $b_1=1$ and $b_2=5$. By
direct computation this implies
$h^1(\Cal O_{\overline X})=0$. Since $C_0 \cdot (C_0+2f) = 0$, the
restriction of $p_*\Cal O_{\overline X}$ to a general member $C$ of
$|C_0+2f|$ is the same as the restriction of $\varphi_*\Cal O_X$ to a
general hyperplane section  of $W$. But then \hypsplit tells that
$b_1=b_2=3$, so we get a contradiction. Therefore the only possibility
left is $a_1=a_2=2$. Then $b_1+b_2=6$, $b_1 \geq 2a_1 -1$ and $b_2
\geq 2a_2 -1$ yields $b_1=b_2=3$. Then by direct computation,
$h^1(\Cal O_{\overline X})=0$, so both $\overline X$ and $X$ are
regular and 3) follows.

\smallskip

 Now we prove 5). Recall that we showed $D_3=C_0$.
We know also that $L_1^{\otimes 2} \otimes \Cal O_Y(-C_0)= \Cal O_Y(D_2)$ and
$L_2^{\otimes 2} \otimes \Cal O_Y(-C_0)= \Cal O_Y(D_1)$. Then, because
of the
values of the
$a_i$s and the
$b_j$s just found, we have that $D_1 \sim D_2 \sim
3C_0+6f$. Recall also that the normality of $\overline X_i$ for
$i=1,2,3$ implies that $D_1+D_2+D_3=D_1+D_2+C_0$ does
not have multiple components. Now we prove the statement in 5),
namely, that
$\overline X$ is the
 normalization of the fiber 
 product over $Y$ of two double covers $p_1^1$ and $p_2^1$ of $Y$, 
  branched along $D_2'=D_2 +
 C_0$ and along $D_1'=D_1 + C_0$ respectively. 
   Let $\hat 
 X @> \hat p >> Y$ be the fiber product over $Y$ of $\overline X_1$ and
 $\overline X_2$.  Let $U=Y - C_0$, $V={\overline X}-F_1-F_2$ and 
 $\hat V={\hat p}^{-1}(U)$. Since $b_1+b_2=b_3$, we have that 
 $(L_1\otimes L_2)|_U=L_3|_U$. Then, by the same reason argued for
 \text{\bipro \ns} (for more details, see \GPsmooth \ns, Proposition 2.7),  
the restriction
 $V @> p|_V >> U$ is a fiber product of the restriction of $p_1^1$ and
 $p_2^1$ to $V$. 
 Thus $V=\hat
 V$. In particular, since $\overline X$ is normal, so is $\hat
V$. 
  Let now $\tilde X$ be the normalization of the reduced part of
 $\hat X$ and $\tilde V$ the open set of $\tilde X$ lying over $U$. Then,
 since $V=\hat
 V$  is normal, $\tilde V=\hat V = V$, so $\tilde X$
 and $\overline X$ are birational. Moreover, $\Cal O_{\tilde X}$ is
 integral over $\Cal O_Y$ and therefore, the integral closure of $\Cal
 O_Y$ in $\Cal K(X)$. Hence $\tilde X=\overline X$.

 \smallskip
 
 Now we prove the converse.   
 Let $Y=\bold F_2$ and let $\overline X @> p >> Y$ be the normalization
 of the fiber product 
 of two double covers of $Y$, one branched along
 $D_1 + C_0$, and the other branched along $D_2 +
 C_0$, $D_1 \sim D_2 \sim 3C_0+6f$ and $D_1 +
 D_2 + C_0$ without multiple components. 
 Let  $\overline{X'} @> p_1 >> Y$ be the double cover of $Y$
 branched along $D_1'=D_1 + C_0$. In fact ${\overline
   X}'=\text{Spec}(\Cal O_Y \oplus \Cal O_Y(-2C_0-3f))$. 
 Then $\overline{X'}$ is normal and locally Gorenstein. 
 Now, since $D_1 +
 D_2 + C_0$ have no multiple components, the double cover $p_2$
 of $\overline{X'}$ branched along 
 $p_1^*D_2$ is normal and is in fact 
$\overline
 X$. 
  We denote $\overline{C_0}=p_1^{-1}C_0$ and, for general $f$,
 $\overline f=p_1^{-1}f=p_1^*f$. Then $p_1^*C_0=2\overline{C_0}$. 
 The canonical bundle of $\overline{X'}$ is 
 $$\omega_{\overline{X'}}=
 p_1^*(\omega_Y \otimes \Cal O_Y(2C_0+3f))=p_1^*\Cal O_Y(-f)=\Cal
 O_{\overline{X'}}(-\overline f)\ .$$ 
 Then ${\overline X}=\text{Spec}(\Cal
 O_{\overline{X'}} \oplus \Cal
 O_{\overline{X'}}(-3\overline{C_0}-3\overline f))$ and  
 the canonical bundle of $\overline X$ is
 $$\displaylines{\omega_{\overline X}=p_2^*(\omega_{\overline{X'}} \otimes \Cal
 O_{\overline{X'}}(3\overline{C_0}+3\overline f))= \cr p_2^*(\Cal
 O_{\overline{X'}}(3\overline{C_0}+2\overline f))=p^*(\Cal
 O_Y(C_0+2f))\otimes p_2^*(\Cal O_{\overline{X'}}(\overline{C_0}) 
 ).}$$
 Recall that $\overline{X'}$ is smooth at every point of 
 $\overline{C_0}$, for $D_1 \cdot C_0=0$. 
 Moreover $p_2$ is \'etale at every point of
 $\overline{C_0}$ and $p_2^*\overline{C_0}=F_1+F_2$, where $F_1$ and
 $F_2$ are two disjoint lines, each of them with self-intersection
 $-1$. Then $\overline X$ is smooth at every point of $F_1$ and $F_2$. 
 
 Let $L=p^*\Cal O_Y(C_0+2f)$. Then $L$ is base-point-free. Using
 projection formula we compare
 $H^0(L)$ and $H^0(\Cal O_Y(C_0+2f))$ and see that they are equal. 
 This means that the morphism induced by $H^0(L)$ factorizes through
 $p$.  On the other hand let $\overline X @> \overline q >> X$ be the
 contraction of $F_1$ and $F_2$. Since $\overline X$ is smooth at every
 point of $F_1$ and $F_2$, and $F_1^2=F_2^2=-1$, $X$ is smooth at the
 images $x_1$ and $x_2$ of $F_1$ and $F_2$. Since $\overline X$ is
 normal with at worst canonical singularities, 
 then so is $X$.  
 We also know that 
$$K_{\overline X}={\overline q}^*K_X + F_1 +
 F_2 \quad \eqno \ncrfour \ .$$
 Since $\omega_{\overline X} =p^*(\Cal
    O_Y(C_0+2f))\otimes p_2^*(\Cal O_{\overline{X'}}(\overline{C_0}))$
    and $F_1+F_2=p_2^*\overline{C_0}$, then ${\overline
      q}^*\omega_X=L$. Moreover $H^0(L)=H^0({\overline q}^*\omega_X)=
    H^0(\overline q_*{\overline 
      q}^*\omega_X)=H^0(\omega_X)$. Then, since $L$ is base-point-free,
    so is 
    $\omega_X$ and the morphism induced by $H^0(L)$ also factorizes
    through the canonical morphism of $X$. Thus we have finally the
    desired commutative diagram:

  $$ \matrix \overline X & @> \overline q >> & X \cr
         @VV p V &  \hskip -.8 truecm @VV\varphi V \cr
 Y & @> q >>& W \cr 
 \endmatrix \ $$
 where $W$ is the cone over a conic inside $\bold P^3$, $q$ is the
 minimal desingularization of $W$, $\varphi$ is the canonical
 morphism of $X$, and, by \ncrfour \ns, $\overline q$ is noncrepant. 
Now the fact that over
$Y-C_0$
 the surface $\overline X$ is a fiber product and \biproconv imply
 that $\Cal K(X) 
 / \Cal K(Y)$ is a Galois extension with Galois group $\bold Z_2 \times
 \bold Z_2$. Now, since $p$ is finite and $\overline X$ is normal, $p$
 is Galois cover with group $\bold Z_2 \times \bold Z_2$ and, by
 \samegrplem \ns,  so is
 $\varphi$. \qed

 \proclaim{\noncrepantcyc} Let $W$ be a
 singular rational normal scroll and let $X 
 @> \varphi >> W$ be a Galois canonical cover with Galois
 group $\bold Z_4$. 
 Let $\overline X$, $Y$, $q$, $\overline q$ and $p$  be as in 
 \defXbar \ns.
 If $\overline q$ is not crepant, then

\item{1)} $W=S(0,2)$; 
\item{2)} $\overline X$ has at worst
 canonical singularities;
\item{3)}  $X$ is regular; 
\item{4)} $\overline X @> \overline q >> X$ is the morphism to the
canonical model of $\overline X$; 
 \item{5)} the morphism 
 $p$ is the composition of two double covers 
 $\overline X_1 @> p_1 >> Y$, branched along a divisor $\Delta_2$, 
 and $\overline X @> p_2 >> \overline X_1$,
 branched along the ramification of $p_1$ and $p_1^*D_1$ and with trace
 zero module $p_1^*\Cal O_Y(-\frac 1 2 (D_1 + C_0) -\frac 1 4 \Delta_2) 
 \otimes \Cal O_{\overline X_1}(\overline C_0)$, where $\overline C_0$ is
 $p_1^{-1}C_0$, and  
 either 
 
 \smallskip
 
 \itemitem{5.1)}  $D_1 \sim C_0+3f$, $\Delta_2 \sim 4C_0+6f$; or

 \itemitem{5.2)} $D_1 \sim 4C_0+9f$, $\Delta_2 \sim 2C_0+2f$

 \smallskip

 \medskip
 
 Conversely, let $\overline X$ be a normal surface with at worst
 canonical singularities and let $Y=\bold F_2$. If $\overline X @> p >>
 Y$ is the composition of two
 double covers $p_1$ and $p_2$ as described in 5) above, 
               then  there exists a commutative diagram like
               \text{\diagramone \ns,}
  where $\varphi$ is the canonical morphism of $X$ and is Galois with
Galois group $\bold Z_4$ and $\overline q$
 is noncrepant.  
 \endproclaim

 \noindent {\it Proof.}  From \singlemma and \thmsinglept \ns, if $G=\bold Z_4$ and $\overline q$ is
 noncrepant, $W=S(0,2)$, so we have 1), and 
$\varphi^{-1}\{w\}$ consists of $1$ or $2$ points. We
 split the argument in two cases accordingly: 
 
 \smallskip
 
 \noindent{\it Case 1:}  
 Cardinality of $\varphi^{-1}\{w\}=2$. In this case, according to
 \text{\singlemma \ns,} $\overline X$ is locally Gorenstein
and  
 $\overline q$ is the blowing up of $X$ at $x_1$ and $x_2$, which are
 smooth points.
Moreover $K_{\overline X}= {\overline q}^*K_X + F_1 + F_2$,
 where $F_1$ and $F_2$ are the exceptional divisors,
 $p^*C_0=2F_1+2F_2$,  and $K_{\overline
   X}$ and $K_X$ are suitable canonical divisors. 
 Then in particular the inertia group of $F_1$ and $F_2$ is $\bold
 Z_2$. 
 On the other hand we have 
 $$K_{\overline X} \equiv {\overline q}^*K_X + F_1 + F_2 \equiv p^*(\frac 3 2
 C_0 + 2f) \quad \eqno \ncrtwo .$$
 Arguing  as in the proof of \thmsinglept \ns,
  we have that
 $$p_*\Cal O_{\overline X}= \Cal O_Y \oplus L_1^* \oplus L_2^* \oplus
 L_3^*$$
 and
  $$\displaylines{L_1 \otimes L_1= L_2 \otimes \Cal O_Y(D_1+D_2) \cr
 L_1 \otimes L_2 = L_3 \otimes \Cal O_Y(D_{2})\cr
 L_1 \otimes L_3 = \Cal O_Y(D_1+D_{2}+D_{3})\cr
 L_2\otimes L_2 = \Cal O_Y(D_{2}+D_{3})\cr
 L_2\otimes L_{3} = L_1 \otimes \Cal O_Y(D_{3})\cr
 L_{3} \otimes L_{3} = L_{2} \otimes \Cal O_Y(D_{1}+D_{3}) \ \cr}$$
 where $D_1$, $D_2$, $D_3$ are effective divisors with neither multiple
 components nor commom components pairwise and $D_2$ is either $0$ or
 $C_0$. Furthermore, by \singleptone (note that in the proof of
 \singleptone we do not use the hypothesis of \thmsinglept that
 $\varphi^{-1}w$ is a single point) we have
 $$2K_{\overline X} \sim p^*(2K_Y + 2L_1 + 2L_2 - D_{2}) \sim p^*(2K_Y +
 2L_3 + D_{2}) \quad \eqno \ncrthree . $$
 Comparing \ncrtwo and \ncrthree 
 we conclude that $D_2 \neq 0$, otherwise 
  $\frac 3 2
 C_0 + 2f$ would be numerically equivalent to a Cartier divisor on
 $Y$. Hence $D_2=C_0$. But this would imply that the
 ramification lying over $C_0$ would have inertia group $\bold Z_4$,
 which contradicts the fact that the inertia group of $F_1$ and $F_2$ is $\bold
 Z_2$. Thus Case 1 does not actually occur, and we have proven that 
there is a unique point $x \in X$ lying over the
vertex of $W$. 
 
 \smallskip
 
 \noindent{\it Case 2.} Cardinality of $\varphi^{-1}\{w\}=1$. In this case
we know
 by \thmsinglept \ns,
 that 
 $\overline X$ is locally $2$-Gorenstein and has at worst rational
 singularities. Recall also that  $e=2$,  $p^*C_0=4F$ with $F$
 isomorphic to $\bold P^1$ and $F^2=-\frac 1 2$. Let $L_1, L_2, L_3,
 D_1, D_2, D_3$ be as in the proof of \thmsinglept \ns, Case 1.2,
  and let $L_i=\Cal
 O_Y(a_iC_0+b_if)$. 
 Then $D_2=C_0$ and $L_1 \otimes L_2 = L_3
 \otimes \Cal O_Y(C_0)$. 
Moreover, 
 $$\displaylines{2K_{\overline X}= {\overline q}^*2K_X + 4F \sim
 {\overline q}^*2K_X + 
  p^*C_0 \cr 
 K_{\overline X} = {\overline q}^*K_X + 2F \equiv {\overline q}^*K_X +
 \frac 1 2 p^* C_0  \cr
 2K_{\overline X} \sim p^*(2K_Y + 2L_3 + C_0) \cr
 K_{\overline X} \equiv p^*(K_Y + L_3 + \frac 1 2 C_0) }$$
 From this, we see that $\omega_Y \otimes L_3 = \Cal O_Y(C_0+2f)$, hence 
 $L_3 = \Cal O_Y(3C_0+6f)$. 
 
 Therefore we have $a_1+a_2=4$ and $b_1+b_2=6$. 
 We examine all possibilities for the $a_i$'s and the $b_j$'s. 
 Recall that $L_1^{\otimes 2} \otimes L_2^*$ and $L_2^{\otimes 2}$ are
 effective, hence   
 $$\displaylines{2a_1-a_2 \geq 0 \cr
 2b_1 -b_2 \geq 0  \quad \rlap \ncrfive \cr
 a_2, b_2  \geq 0}$$
 
 Then $a_1, b_1 \geq 0$ also.
 On the other hand, 
  $D_2+D_3$ has no multiple components, so in particular, 
 the components of the fixed part of $|L_2^{\otimes 2}|$ have
 multiplicity $1$, and hence $C_0$ appears with at most
 multiplicity $1$ in  the fixed part of $|L_2^{\otimes 2}|$. 
 Likewise $D_{1}$
 and $D_{2}$ do not have common components, and in particular $C_0$ is
 not a common component of $D_{1}$
 and $D_{2}$. Since $D_{2}=C_0$, we have that $|L_1^{\otimes 2}
 \otimes L_2^* \otimes \Cal O_Y(-C_0)|$ does not have $C_0$ as fixed
 component.  Since $D_{3}$ does not contain $C_0$ either, $|L_2 \otimes
 L_3 \otimes L_1^*|$ does not have $C_0$ as fixed component. 
 This yields the inequalities 
 $$\displaylines{b_2 \geq 2a_2 -1 \cr
 2b_1 - b_2 \geq 2(2a_1 -a_2 -1) -1 \ .\quad \ncrsix \cr 
 }
 $$

 We start ruling out possible values for $a_1$. 
 The value $a_1=0$ is not possible for in that case $a_2=0$, but $a_1+a_2=4$. 
 Suppose now that $a_1=1$. Then, by \ncrfive \ns,  $0 \leq a_2 \leq 2$,
hence
$a_1+a_2
 \leq 3$ and we reach
 again a contradiction. Now, if $a_1=4$, then $a_2=0$ and  $2b_1 -b_2
 \geq 13$.  This is not
 possible, since $b_2 \geq 0$ and $b_1+b_2=6$.
 
 Therefore the only values for $a_1$ which are  still possible are $2$
and $3$. If 
 $a_1=2$,  then $a_2=2$ and  by \ncrsix \ns, $b_2 \geq 3$  and hence $b_1
\leq
 3$. Then, again by \ncrsix \ns, we should have $b_1=b_2=3$. 
 This corresponds to 5.1) in the statement.  
 
 Finally, if $a_1=3$, then $a_2=1$, $b_2 \geq 1$ and $b_1 \leq 5$. 
 But $2b_1-b_2
 \geq 7$ by \ncrsix \ns, hence we should have  $b_1=5$, $b_2=1$. This
corresponds to
 5.2) in the statement. 
 
 \smallskip
 
 Now we finish the description of $p$ separately for $a_1=2$ and
 $a_1=3$. In case $a_1=2$, recall that $a_1=a_2=2$ and $b_1=b_2=3$
 and  
 we have 
 $$p_*\Cal O_X=\Cal O_Y \oplus \Cal O_Y(-2C_0-3f) \oplus \Cal
 O_Y(-2C_0-3f) \oplus \Cal O_Y(-3C_0-6f)$$
 and $D_{1} \sim C_0+3f$, $D_{2}=C_0$ and $D_{3} \sim 3C_0+6f$.

 Recall 
 also that $C_0$ is contained neither in $D_{1}$ nor in $D_{3}$.  
 Let $U=Y-C_0$ and $V=\overline X - F$. Abusing the notation,
 we will call also $D_{1}$ and $D_{3}$ the restrictions of $D_{1}$
 and $D_{3}$ to $U$. 
  Likewise we
 call $L_1, L_2, L_3$ to the restrictions of $L_1, L_2, L_3$ to $U$.  
 
 Then, on $U$ we have:  
 $$\displaylines{L_1 \otimes L_1 = L_2 \otimes \Cal O_Y(D_{1}) \cr
 L_1 \otimes L_2 = L_3 \cr
 L_1 \otimes L_3 = \Cal O_Y(D_{1}+D_{3})\cr
 L_2\otimes L_2 = \Cal O_Y(D_{3})\cr
 L_2\otimes L_3 = L_1 \otimes \Cal O_Y(D_{3})\cr
 L_3 \otimes L_3 = L_2 \otimes \Cal O_Y(D_{1}+D_{3})\cr}$$
 and $p|_V$ is the composition of $U'@> \pi_1 >> U$, a double cover of
$U$
 branched along $D_{3}$ and $\pi_2$, a double cover of $U'$ branched
 along $\pi_1^*D_{1}$ 
 and the ramification of $\pi_1$.
 
 On the other hand we consider the double cover $\overline X' @> p_1
 >> Y$ branched along $C_0 + D_{3}$ and the cover $\hat X @> p_2
 >> \overline X'$, branched along the ramification of $p_1$ and
 $p_1^*D_{1}$ and with trace zero module $p_1^*\Cal O_Y(-C_0-3f)
 \otimes \Cal O_{\overline X_1}(-\overline C_0)$.
 Let $\hat X_{\text{norm}}$ be the normalization of
 $\hat X$. First note that the open set of $\hat X$ lying over $U$ is
 equal to $V$, which is normal since $\overline X$ is. 
On the points of $\hat X$ lying over $C_0$ we see only one
 singularity of type $A_1$: the point lying over the intersection of
 $D_{1}$ and $C_0$. Indeed,  recall that $C_0 \cdot D_{1}=1$, hence the
 intersection is transversal and so is the intersection of
 $p_1^*D_{1}$ and the ramification of $p_1$ lying over $C_0$. 
 Hence $\hat X$ is normal everywhere. By construction, the open set
of $\hat X$ lying over $U$ is $V$, so $\overline X$ and $\hat X$ are
birational. Since $\hat X$ is normal and integral over $Y$, $\hat X$ is
in fact the integral closure of $\Cal O_Y$ in $\Cal K(X)$, so in fact
$\overline X=\hat X$. 
Thus we have seen that  
 $p$ is the composition of two double covers 
 $\overline X_1 @> p_1 >> Y$ branched along a divisor $\Delta_2 =D_3 +
 C_0 \sim 4C_0+6f$ 
 and $\overline X @> p_2 >> \overline X_1$,
 branched along the ramification of $p_1$ and $p_1^*D_1$, where   
 $D_1 \sim C_0+3f$. We have also seen
  that the trace zero module of $p_2$  is $p_1^*\Cal O_Y(-C_0-3f)
  \otimes \Cal O_{\overline X_1}(-\overline C_0)$. This proves 5.1). 
 Now we prove 2). We know that outside $F$, the surface $\overline X$ and $X$
 are isomorphic so, outside $F$, $\overline X$ has canonical
 singularities by hypothesis. On the other hand, we have seen that
 $\overline X=\hat X$, and that the points of $F$ are smooth points of
 $\hat X$ except for one point which is an $A_1$ singularity, which is
 a canonical singularity. Thus
2) is proven in case 5.1).

Now, since $p$ is the composition of two double covers, and since we
know its trace zero modules we can easily see that
$$\omega_{\overline X} = p^*\Cal O_Y(C_0+2f) \otimes p_2^*\Cal
O_{\overline X_1}(\overline C_0) \ .\quad \eqno \ncrseven $$
Then $K_{\overline X}
\cdot F=-1$ and $K_{\overline X}$ intersects strictly
positively every other curve of $\overline X$, so $X$ is not only 
minimal but it is also the canonical model  of $\overline X$.
  This ends
the proof of  4) in case 5.1). 
 
 \smallskip
 
 Finally, we describe $p$ if $a_1=3$. Then $a_1=3, a_2=1, b_1=5, b_2=1$
 so  
 we have
 $$p_*\Cal O_X=\Cal O_Y \oplus \Cal O_Y(-3C_0-5f)  \oplus \Cal
 O_Y(-C_0-f) \oplus \Cal O_Y(-3C_0-6f) \ .$$
 In this case, $D_{2}=C_0$, $D_{1} \sim 4C_0+9f$, $D_{3} \sim
 C_0+2f$. Like before $C_0$ is contained in neither $D_{1}$ nor
 $D_{3}$. Then we can argue as before to show that $p$ is the
 composition of a double cover $\overline X' @> p_1 >> Y$ branched
 along $\Delta_2=C_0 + D_{3} \sim 2C_0+2f$ and a double cover
 $\overline X @> p_2 >> 
 \overline {X'}$, branched along the ramification of $p_1$ and
 $p_1^*D_{1}$ and with trace zero module $p_1^*\Cal O_Y(-2C_0-5f)
 \otimes \Cal O_{\overline X_1}(-\overline C_0)$. This proves the
description in 5.2).
 Again $D_{1} \cdot C_0=1$, and $D_{3} \cdot
 C_0=0$, hence there is only one point on $F$, the one lying
 over $D_{1} \cap C_0$, which is singular, and its singularity is of
 type $A_1$. Then, by the same reason as before, $\overline X$ is
 Gorenstein and therefore has canonical singularities, so we have 2)
in case 5.2). The proof of 4) in this case is as well as in the case 5.1). 

 \smallskip
 
 To prove that $\overline X$ is regular 
 we need just to use
 the splitting of $p_*\Cal O_{\overline X}$, which is determined by the
values of the $a_i$s and the $b_j$s corresponding to 5.1) and 5.2), and
compute the cohomology. Since $\overline X$ is birational to $X$ 
and both have rational singularities, $X$ is also regular and hence
3) is proven.

 \smallskip

 Now we prove the converse. Let $\overline X$ be a normal surface with
 at worst canonical singularities and let $\overline X @> p >> Y$ be
the composition of two double covers  
 $\overline X_1 @> p_1 >>Y$ branched along a divisor $\Delta_2$ 
  and $\overline X @> p_2 >> \overline X_1$,
 branched along the ramification of $p_1$ and $p_1^*D_1$,
 where
 $D_1$ and $\Delta_2$ satisfy condition 5.1) or 5.2) of the statement.
Let $\overline C_0$ be the inverse image of $C_0$ by $p_1$ and let $F$ be
 the inverse image of $C_0$ by $p$. Let finally 
 $p_1^*(\Cal O_Y(C_0+3f)) \otimes \Cal O_{\overline X_1}(\overline
 C_0)$ or $p_1^*(\Cal O_Y(2C_0+5f)) \otimes \Cal O_{\overline X_1}(\overline
 C_0)$ be the trace zero module of $p_2$ accordingly. Then one easily
obtains as before the formula \ncrseven for the canonical of
$\overline X$.
Thus $L=p^*(\Cal
O_Y(C_0+2f))$ is the
 free part of $\omega_{\overline X}$. 
 
 We compare now $H^0(L)$ and $H^0(\Cal
 O_Y(C_0+2f))$. In the first place, using projection formula we see
 $H^0(p_1^*\Cal 
 O_Y(C_0+2f))=H^0(\Cal
 O_Y(C_0+2f))$. Similarly $H^0(p^*\Cal 
 O_Y(C_0+2f))=H^0(p_1^*\Cal 
 O_Y(C_0+2f))$. Thus the morphism induced by $|L|$ factorizes as $q
 \circ p$, where $Y @> q >> S(0,2)$ is the morphism induced by
 $|C_0+2f|$.

Now, since $C_0$ is a component of $\Delta_2$, by construction of $p$, 
$p^*C_0=4F$, where $F$ is a smooth line. By Stein factorization, $q \circ
p$ factorizes as the composition of $\overline X @> \overline q >> X$
followed by
$X @>
\varphi
>> W$, where, by the commutativity $q \circ p = \varphi \circ \overline
q$, $\overline q$ contracts only
$F$ and $\varphi$ is finite. Since by hypothesis $\overline X$ is
normal and $\Delta_2 \cdot (\Delta_2 - C_0)=0$, $\Delta_2$ is
smooth along $C_0$. Since in addition $D_1 \cdot C_0=1$, $F$ contains
only one singular point of $\overline X$, which is of type $A_1$.
Contracting $F$ give raise
to a smooth point in $X$. Then, since $\overline X$ has canonical
singularities so does $X$. 
Now, since $p^*C_0=4F$, we have 
$F^2=-1/2$. By \ncrseven we have also $\omega_{\overline X}
\cdot F=-1$, so 
$$K_{\overline X}=\overline q^*K_X + 2F$$ for suitable canonical
divisors, so $\overline q$ is noncrepant.  Comparing this formula with
\ncrseven yields $\overline q^*\omega_X=p^*\Cal O_Y(C_0+2f)$, so
$H^0(L)=
H^0(\overline q^*\omega_X)=
H^0(\omega_X)$, so $\omega_X$ is base-point-free, for so is $L$, and
in the factorization $q \circ p = \varphi \circ \overline q$, $\varphi$ is
in fact the canonical morphism of $X$.

Finally by \Galcycliconv 
 $p$, and therefore $\varphi$, are Galois covers with Galois group
$\bold Z_4$.
 \qed

\medskip 

We finish this section summarizing the splitting of $p_*\Cal
O_{\overline X}$ for all the surfaces $\overline X$ which appear in
\classing \ns.
Even though we already showed there that $\overline X$ is regular,
 the reader can check this fact at once by looking at
the next corollary:

 \proclaim{\singsplit} Let $W=S(0,2)$, let $X @> \varphi >> W$ be
a
 quadruple Galois canonical 
 cover and let $\overline X$ and $p$ be as in \defXbar \ns. Then
the vector bundle $p_*\Cal O_{\overline X}$ splits as follows:

\item{1)} $p_*\Cal O_{\overline X} = \Cal O_Y \oplus \Cal
O_Y(-C_0-3f)
\oplus
\Cal O_Y(-2C_0-3f)  \oplus \Cal O_Y(-3C_0-6f)$ if $X$ is as in
\crepantbid or \crepantcyc \ns, i.e., if $\overline q$ is crepant. 

\item{2)}       $p_*\Cal O_{\overline X} = \Cal O_Y \oplus \Cal
O_Y(-2C_0-3f)
\oplus
\Cal O_Y(-2C_0-3f)  \oplus \Cal O_Y(-3C_0-6f)$, if $X$ is as in
\noncrepantbid or in  \noncrepantcyc \ns, 5.1).

\item{3)} $p_*\Cal O_{\overline X} = \Cal O_Y \oplus \Cal
O_Y(-C_0-f)
\oplus
\Cal O_Y(-3C_0-5f)  \oplus \Cal O_Y(-3C_0-6f)$, if $X$ is as in
\noncrepantcyc \ns, 5.2).

 \endproclaim

{\it Proof.} The corollary follows from  the values for $a_i$s and $b_j$s
found in the proofs of \classing \ns.  \qed

\heading 5. Singularities of quadruple Galois canonical covers and
  examples. \endheading

In this section we describe further the surfaces
$X$ and $\overline X$ classified in
\classing and the morphism $\overline q$. We focus especially 
in the study of the singularities of $X$ and $\overline X$.

\proclaim{\numbsing } Let $W=S(0,2)$ and let $w$ be its vertex. Let
$X @>
\varphi >> W$ be a quadruple Galois canonical 
 cover with Galois group $G$ and let $\overline X$ and $\overline q$
  be as in 
\defXbar
\ns.

\item{1)} If $G=\bold Z_2 \times \bold Z_2$ and $\overline q$ is
crepant (i.e., $\varphi$ is as in \crepantbid \ns), then there is only one
point $x$ lying over $w$. Moreover, $x$ is at best an  $A_1$ 
singularity and in general an $A_l$ singularity. Moreover
 $\overline X
@>
\overline q >> X$ is the minimal desingularization of $x$ if $x$ is of
type $A_1$ and a partial desingularization of $x$ (that consists of two
consecutive blowing ups) otherwise.

\item{2)} If $G=\bold Z_2 \times \bold Z_2$ and $\overline q$ is
noncrepant (i.e., $\varphi$ is as \noncrepantbid \ns), then there are
only two points $x_1$ and $x_2$ lying over $w$, they are smooth and
$\overline X @> \overline q >> X$ is the blowing up of $X$ at $x_1$ and
$x_2$. 

\item{3)} If $G=\bold Z_4$ and $\overline q$ is
crepant (i.e., $\varphi$ is as \crepantcyc \ns), then 

\itemitem{3.1)} there is only one point $x$ lying over $w$. Moreover, $x$ is
a
$D_4$ singularity and $\overline X @> \overline q >> X$ is the blowing
up of
$X$ at
$x$;

\itemitem{3.2)} $X-\{x\}$ is singular and the mildest possible set of
singularities of $X-\{x\}$
 consists of $9$ $A_1$ singularities. 

\item{4)} If $G=\bold Z_4$ and $\overline q$ is
noncrepant (i.e., $\varphi$ is as \noncrepantcyc \ns), then 

\itemitem{4.1)} there is only one point $x$ lying over $w$ and $x$ is
smooth.

\itemitem{4.2)}  The morphism $\overline X @> \overline q >> X$ is the
contraction of a smooth line $F$, which is $p^{-1}C_0$. The line $F$
consists of smooth points of $\overline X$ and an $A_1$
singularity and its self-intersection is $F^2=-\frac
1 2$. Moreover $K_{\overline
X}={\overline q}^*K_X + 2F$;

\itemitem{4.3)} $X-\{x\}$ is singular and the mildest possible set of
singularities of 
$X-\{x\}$ consists of $9$ $A_1$ singularities. 

\endproclaim

{\it Proof.} First we prove 1). We use the description of the branch
divisors of
$p$ given in
\crepantbid \ns, 5). Since $D_1 \sim 4C_0+6f$, $C_0$ is a
component of $D_1$, so we can write $D_1=D_1'+C_0$. Since
$\overline X$ is normal, $C_0$, $D_1'$ and $D_2$ have no common
components. Since  $C_0 \cdot (3C_0+6f)=0$, $D_1$ does not meet
$C_0$ so $\overline X_2$ is smooth at the points (which are in the
ramification locus of $p_2$) lying over $C_0$. On the other hand 
$D_2
\cdot C_0=2$. Then $D_2$ can meet $C_0$ transversally or not. In
the second case $D_2$ can be smooth at the intersection point
with $C_0$ or can have an $A_k$ singularity. All this means that the
inverse image of
$C_0$ by
$p$ is a $-2$ cycle $Z$ consisting of one smooth $-2$-line or two lines
meeting at one point $\overline x$, which is either a smooth point of
$\overline X$ or an $A_k$ singularity, and the points of $Z$ are smooth
points of $\overline X$ except maybe $\overline x$. In any case
contracting $Z$, as $\overline q$ does, gives rise to a unique point $x$
lying over $w$, which is a singularity of type $A_{k+2}$ if $\overline
x$ is of type $A_k$, is of type $A_2$ if $\overline x$ is smooth but
$D_2$ and $C_0$ do not meet transversally and is of type $A_1$ if 
$D_2$ and
$C_0$ meet  transversally.
This can be easily seen resolving
$\overline x$ if necessary and looking at the total transform of $Z$.

\smallskip 

Part 2) was already proven in \singlemma so we now prove 3.1). As
argued before, since $\overline X$ is normal, $D_2=D_2' + C_0$ and
$D_2$ does not meet $C_0$. On the other hand $D_1 \sim 3f$, and
since $\overline X$ is normal $D_1$ consists of $3$ distinct fibers of
$\bold F_2$. Thus, if we call $E=p_1^{-1}C_0$, all points of $E$ are
smooth points of $\overline X_1$ and, near $E$ $p_2$ is branched
along $E+p_1^*D_1$. Such curve is smooth at all points of $E$ except at
$3$ distinct points which are $A_1$ singularities. Then if we call
$F=p^{-1}C_0$, $F$  is a smooth line
 with $F^2=-1/2$ lying over $C_0$ and all the points of $F$ are smooth
points of $\overline X$ except three distinct points $\overline
x_1$, $\overline x_2$ and $\overline x_3$ which are
$A_1$ singularities. Now $\overline X @> \overline q >> X$ contracts
only
$F$ and therefore gives rise to a 
single point $x$ lying over $w$, and $x$ is a singularity of type $D_4$.
The last claim is inmediate once we resolve $\overline X$ at $\overline
x_1$, $\overline x_2$ and $\overline x_3$, since the total transform of
$F$ is the $-2$-cycle which appears in the minimal desingularization
of a $D_4$ singularity. In fact, $\overline q$ is a partial desingularization
of $x$ consisting in blowing up $X$ at $x$ once. 

\smallskip

Now we prove 4.1). The argument is similar to 3.1). The
description of the branch divisors of $p$ given in \noncrepantcyc
\ns, 5.1) and 5.2) and the fact that $\overline X$ is normal implies
that $\Delta_2=D_3+C_0$ and $D_1 + D_3 + C_0$ does not have
multiple components. Since $\Delta_2 \cdot C_0=0$, $\overline X_1$ is
smooth along $\overline C_0=p_1^{-1}C_0$. 
Since $D_1 \cdot C_0=1$, in both 5.1) and
5.2), then $D_1$ and $C_0$ meet transversally at a point. Since, near
$\overline C_0$, 
$p_2$ is branched at $\overline C_0+p_1^*D_1$,
there is only one singular point $\overline x$ of $\overline X$ lying on
$F=p^{-1}C_0$, and $\overline x$ is an $A_1$ singularity. On the
other hand, $F$ is a smooth line with $F^2=-\frac 1 2$ as in 3.1)
Again we resolve $\overline x$ and the total trasform $T$ is a cycle
with self-intersection $-1$ consisting of two smooth lines meeting
transversally at one point and with self-intersections $-1$ and
$-2$. Then the contraction of $T$ is a smooth point, and so
$\overline X @> \overline q >> X$ contracts $F$ to a unique point $x$
lying over $w$, and $x$ is smooth. Finally $K_{\overline X}=\overline
q^* K_X + 2F$ follows from \ncrseven \ns, since $p^*\Cal
O_Y(C_0+2f)=\overline q^*K_X$ and $p_2^*\overline C_0=2F$. 

\smallskip

Finally, we prove 3.2) and 4.2). 
The mildest singularities in
$\overline X - F$ occur when $D_2 - C_0$ and $D_1$ in \crepantcyc
\ns, 5) and $\Delta_2-C_0$ and $D_1$ in \noncrepantcyc \ns, 5.1) and
5.2) meet transversally. Since $(D_2-C_0) \cdot D_1= (\Delta_2-C_0)
\cdot D_1= 9$, if the intersection is transversal $\overline X -F$ has
$9$ singular points which are of type $A_1$ and so has $X - \{x\}$. 
\qed

\smallskip

 We end this section showing the existence of surfaces $X$ like those
 classified in \classing \ns. We use the notation of \numbsing \ns.

 \proclaim{\excrepantbid} There exist families of quadruple canonical 
covers $X @> \varphi >> W$ as in 
 \crepantbid with $X-\{x\}$ smooth and $x$
an
$A_1$ singularity. 
 \endproclaim
 
 {\it Proof:} These families were constructed in \GPtrans \ns, Example
 3.7.
 \qed

 \proclaim{\excrepantcyc} There exist families of quadruple canonical 
 covers as in 
 \crepantcyc and \noncrepantcyc with $X-\{x\}$ smooth except for
$9$
$A_1$ singularies
 \endproclaim
 
 {\it Proof.} According to the converse part of \crepantcyc and
\noncrepantcyc we just have
 to construct two double covers $\overline X @> p_2 >> \overline X_1$
 and $\overline X_1 @> p_1 >> Y$ ($Y=\bold F_2$) branched along
 suitable divisors. For surfaces as in \crepantcyc we choose $D_1 \sim
3f$ consisting of three distinct fibers and we choose
 $D_2=D_2' + C_0$ with $D_2' \sim 3C_0+6f$ and so that
$D_1+D_2'+C_0$ has no multiple components. Then the  cover
$\overline X @> p >> Y$ constructed according to \crepantbid is
normal. 
We choose furthermore $D_1$ and $D_2'$ so that $D_2'$ is smooth and
$D_1$ and
$D_2'$ meet transversally. All this can be achieved because $D_2' \sim
3C_0+6f$ is
 base-point-free.  Then as argued in \numbsing $\overline X - F$ has
only
$9$
 singularities, which are of type $A_1$, and has therefore canonical
 singularities, and so has $X$ according to \crepantcyc \ns. Note that
 these are the examples we would obtain in \exampcyc if we allowed
 $m=e=2$. Allowing $D_1+D_2$ to have worse singularities one can
construct $X$ with worse singularities.

\smallskip

To construct surfaces as in \noncrepantcyc again we construct $p_1$
and
$p_2$ following the guidelines in
 the converse part of \noncrepantcyc \ns, 
choosing $D_1$ and $\Delta_2$
so that $D_1+D_3+C_0$ has no multiple components and so
that $D_1$ and $D_3$ are smooth and meet tranversally. This can be
achieved in both cases 5.1) and 5.2) of \noncrepantcyc \ns, 
because $C_0+3f$, $3C_0+6f$, $4C_0+9f$ and $C_0+2f$ are
base-point-free.
 This assures us
 that $\overline X-F$ has only $9$ singular points, which are of type
 $A_1$, so in particular $\overline X$ has canonical singularities,
 and, according to \noncrepantcyc \ns, so does $X$. Allowing $D_1 +
 \Delta_2$ to have worse singularities, one can construct $X$ with
 worse singularities. 
\qed

 \proclaim{\exampnoncrepbid} There are families of smooth quadruple canonical 
covers $X @> \varphi >> W$ as in  
 \noncrepantbid \ns.
 \endproclaim

 {\it Proof.} Let $Y=\bold F_2$. 
 By the converse part in \noncrepantbid we construct
 $\overline X$  as the
 normalization of the fiber product of two double covers 
 $\overline X_1 @> p_1 >>
 Y$ and $\overline X_2 @> p_2 >> Y$, branched respectively along
 divisors $D_1'$ and $D_2'$, where $D_i'=D_i+C_0$ and 
 $D_1 \sim D_2 \sim 3C_0+6f$. We also take $D_1+D_2+C_0$ without
 multiple components. This can be done because $3C_0+6f$ is
 base-point-free. We also take $D_1$ and $D_2$ smooth
and intersecting
 transversally, which again can be achieved because $3C_0+6f$ is
 base-point-free. From the description made in the proof of
 \noncrepantbid of the normalization of the
 fiber product as a composition of two double covers, the second cover
 is 
 branched along the pullback of $D_2$ so we see that $\overline X$ is
 smooth. We also know by \noncrepantbid that the canonical 
morphism of
 $\overline X$ only contracts the inverse image of $C_0$ to two
smooth points in $X$, so  $X$ is a smooth
 surface. If we allow worse singularities in $D_1 + D_2$, then
 examples of singular $X$ can be constructed. 
 \qed

\heading References \endheading

 \item{\Ca} F. Catanese, On the moduli spaces of surfaces of general type. 
 J. Differential Geom. {\bf 19} (1984), no. 2, 483--515.
  
 \item{\CKM} H. Clemens, J. Koll\'ar and S. Mori, {\it Higher dimensional
 complex geometry}, Asterisque {\bf 166} (1988).
 
 \item{\GPtrans} F.J. Gallego and B.P. Purnaprajna, {\it On the
     canonical ring 
 of covers of surfaces of minimal degree}, 
 Trans. Amer. Math. Soc. {\bf 355} (2003), 2715-2732.

\item{\GPsmooth} F.J. Gallego and B.P. Purnaprajna, {\it
    Classification of quadruple Galois canonical covers, I}, preprint.

\item{\CR } F.J. Gallego and B.P. Purnaprajna, {\it Classification of
    quadruple canonical covers: Galois case}, 
to appear in C. R. Math. Acad. Sci. Soc. R. Can.

 \item{\HM} D. Hahn and R. Miranda, {\it Quadruple covers of algebraic
 varieties}, J. Algebraic Geom. {\bf 8} (1999), 1--30.
 
 \item{\Hoone} E. Horikawa, {\it Algebraic surfaces of general type with
 small $c^2_1$, I}, Ann. of Math. (2) {\bf 104} (1976), 357--387.

 \item{\Hothree} E. Horikawa, {\it Algebraic surfaces of general type with
 small $c\sp{2}\sb{1}$, III}, Invent. Math. {\bf 47} (1978), 209--248.

 \item{\Kon} K. Konno, {\it Algebraic surfaces of general type with $c_1^2
 =3p_g-6$}, Math. Ann. {\bf 290} (1991), 77--107.

\item{\Pa} R. Pardini, {\it Abelian covers of algebraic varieties},
  J. Reine Angew. Math.  {\bf 417}  (1991), 191--213.

 \item{\Pu} B.P. Purnaprajna, {\it Geometry of canonical covers with 
 applications to Calabi-Yau threefolds}, in Vector Bundles and 
 Representation Theory, Eds. D. Cutkosky et al., 
Contemporary Mathematics Series AMS {\bf
   322} 
(2003), 107--124.

 \enddocument